\documentclass[leqno,draft,12pt]{amsart}
\usepackage{amssymb}
\usepackage{amsmath}
\usepackage{enumerate}
\usepackage{amsfonts}
\newtheorem{theorem}[equation]{Theorem}
\newtheorem{corollary}[equation]{Corollary}
\newtheorem{lemma}[equation]{Lemma}
\newtheorem{proposition}[equation]{Proposition}
\newtheorem{remark}[equation]{Remark}

\numberwithin{equation}{section} \DeclareMathOperator{\tr}{Tr}
 
\newcommand{\e}{\text{\bf E}}
\newcommand{\p}{\text{\bf P}}
\newcommand{\N}{\mathbb{N}}
\newcommand{\Z}{\mathbb{Z}}
\newcommand{\R}{\mathbb{R}}
\newcommand{\Q}{\mathbb{Q}}
\newcommand{\C}{\mathbb{C}}

\newcommand{\pr}{\mathbb{P}}
\newcommand{\T}{\mathbb{T}}
\newcommand{\m}{M_{inv}(d,\Z)}

\begin{document}

\title[Semigroup actions on tori]{Semigroup actions on
tori and stationary measures on projective spaces}
\subjclass[2000]{54H20, 22E40, 60J05, 60B15} \keywords{Limit set,
proximal and expanding element, toral automorphism, random walk,
projective space, stationary measure}
\begin{abstract}
Let $\Gamma$ be a sub-semigroup of $G=GL(d,\mathbb R),$ $d>1.$ We
assume that the action of $\Gamma$ on $\R^d$ is strongly
irreducible and that $\Gamma$ contains a proximal and expanding
element. We describe contraction properties of the dynamics of
$\Gamma$ on $\R^d$ at infinity. This amounts to the consideration
of the action of $\Gamma$ on some compact homogeneous spaces of
$G,$ which are extensions of the projective space $\pr^{d-1}.$ In
the case where $\Gamma$ is a sub-semigroup of $GL(d,\R)\cap
M(d,\Z)$ and $\Gamma$ has the above properties, we deduce that the
$\Gamma$-orbits on $\T^d=\R^d\slash\Z^d$ are finite or dense.
\end{abstract}

\author[Y. Guivarc'h]{Yves Guivarc'h}
\address{IRMAR, Universit\'e de Rennes 1\\
Campus de Beaulieu\\
35042 Rennes cedex, France} \email{yves.guivarch@univ-rennes1.fr}
\author[R. Urban]{Roman Urban}
\address{Institute of Mathematics\\
Wroclaw University\\
Plac Grunwaldzki 2/4\\
50-384 Wroclaw, Poland} \email{urban@math.uni.wroc.pl}
\thanks{The second author was partially supported by RTN Harmonic Analysis
and Related Problems contract HPRN-CT-2001-00273-HARP and in part
by KBN grant 1PO3A01826.} \dedicatory{Dedicated with admiration to
Hillel Furstenberg on the occasion of his 70th birthday}
\maketitle

\section{Introduction and main results}\label{introduction}
Let $\Gamma$ be a multiplicative semigroup of integers. The
semigroup $\Gamma$ is said to be lacunary if the members
$\{\gamma\in\Gamma:\gamma>0\}$ are of the form $\gamma_0^k,$
$k\in\N,$ $\gamma_0\in\N^*.$ Otherwise $\Gamma$ is non-lacunary.
In 1967 in \cite{F} Furstenberg proved that if $\Gamma$  is a
non-lacunary semigroup of integers and $\alpha$ is an irrational
number, then the orbit $\Gamma\alpha$ is dense modulo 1. The
problem of approximating a number $\theta$ modulo 1 by numbers of
the form $q\alpha,$ where $\alpha$ is a fixed irrational and $q$
varies in a specified subset $Q\subset\N$ was considered by Hardy
and Littlewood in \cite{HL} for various subsets $Q$ of $\N.$ In
particular, the result of Furstenberg above can be considered as a
generalization of a theorem of \cite{HL}, which asserts that if
$r$ is a positive integer and $\alpha$ is an irrational number,
the set $\{q^r\alpha:q\in\N\}$ is dense modulo 1; furthermore,
this result draws attention to the role of the multiplicative
structure of $Q$ in Diophantine approximation, hence of the role
of the corresponding dynamical properties of endomorphisms of
$\T=\R\slash\Z.$ Hence, one is led, more generally, to consider
separately the properties depending on the multiplicative
structure of $\Gamma$ and the properties depending on the additive
structure implied in reduction modulo 1. In this direction, a
generalization of Furstenberg's result to a commutative semigroup
$\Gamma\subset\m:=GL(d,\R)\cap M(d,\Z),$ where $M(d,\Z)$ is the
set of $d\times d$ matrices with integer entries, acting by
endomorphisms on the torus $\T^d=\R^d\slash\Z^d $ was given by
Berend in \cite{B}.

Following \cite{B}, we say that the semigroup of endomorphisms of
a $d$-di\-men\-sio\-nal torus $\mathbb T^d$ has ID-property (cf.
\cite{B, 03a,03b}) if the only infinite closed $\Gamma$-invariant
subset of $\mathbb T^d$ is $\mathbb T^d$ itself. (ID-property
stands for \textit{infinite invariant is dense}.) Berend in
\cite{B} gave necessary and sufficient conditions in arithmetical
terms for a commutative semigroup to have ID-property.

On the other hand, starting from \cite{B} and  \cite{F}, Margulis
in \cite{Ma} asked for necessary and sufficient conditions on
sub-se\-mi\-group $\Gamma\subset \m$ in order that the
$\Gamma$-orbit closures on $\T^d$ are finite union of manifolds.
We observe that it follows from general results of Dani and
Raghavan on linear actions  \cite{DR} that the orbits of
$\Gamma=SL(d,\Z)$ acting on $\T^d$ are finite or dense. In this
direction a detailed study of $\Gamma$-orbits in $\R^d$ of a
general subgroup $\Gamma\subset SL(d,\R)$ was developed by Conze
and Guivarc'h in \cite{CG01}. The homogeneity at infinity of
$\Gamma$-orbits was pointed out there as well as the role of
"$\Gamma$-irrational" vectors in the construction of limit points
of $\Gamma$-orbits, if $\Gamma$ is a general subgroup of
$SL(d,\R).$

Some results in direction of general question of Margulis has been
obtained recently. Muchnik in \cite{03a} proved that if the
semigroup $\Gamma$ of $SL(d,\mathbb Z)$ is Zariski dense in
$SL(d,\mathbb R)$, then $\Gamma$ acting on $\T^d$ has ID-property.
Starkov in \cite{S1} proved the same result in case $\Gamma$ is a
strongly irreducible subgroup of $SL(d,\Z).$ In the next paper
\cite{03b} Muchnik generalized the results of Berend to semigroups
of $\m.$ In the same time Guivarc'h and Starkov in \cite{GS}
derived an important part of Muchnik's result using different
methods, based on \cite{CG00,CG01}. We observe that in \cite{GS},
the property $\Gamma\subset SL(d,\Z)$ is used only when additive
aspects connected with reduction modulo one comes into the play.
It turned out that the property of $\Gamma$-orbits in $\R^d$ which
is responsible for density in $\T^d$ is "thickness" at infinity of
$\Gamma$-orbits (see Theorem \ref{ot} and the comment). Hence,
this property can be studied separately in full generality;
$\Gamma$ is then a general sub-semigroup of $GL(d,\R)$ and the use
of boundaries and random walks is natural in this context.

In this paper we consider this problem in a simplified setting, we
give a self-contained exposition of some of the methods developed
in \cite{CG00,CG01,GS} in the more general context of random walks
and linear actions, and we use the results to prove ID-property in
our setting. We prove also some new results for actions on tori
and on certain compact $G$-factor spaces of $\R^d.$

The general idea is to lift the automorphisms of the torus $\T^d$
to its universal cover $\R^d$ and to study the action of the lifts
at infinity. The action of $\Gamma$ at infinity can be expressed
in terms of some compact homogeneous spaces of $GL(d,\R)$ which
are closely related to the projective spaces $\pr^{d-1}.$ The
random walk framework allows to take into account the global
semigroup asymptotic behavior in terms of stationary measures and
convergence to them. As in \cite{F1} and \cite{G} the results can
be used to obtain topological properties of the $\Gamma$-action.
Furthermore, this general framework allows us to obtain a series
of informations on linear actions which are of interest in other
problems.

Before we state the results we need to introduce some notions. A
matrix $\gamma\in GL(d,\R)$ is said to be \textit{proximal} if it
has an eigenvalue $\lambda_\gamma$ such that
$|\lambda_\gamma|>|\lambda|$ for all other eigenvalues $\lambda$
of $\gamma.$ A matrix $\gamma$ is said to be \textit{expanding} if
it has an eigenvalue $\lambda$ such that $|\lambda|>1.$

Let $\Gamma$ be a sub-semigroup of $GL(d,\R).$ The $\Gamma$-action
on $\R^d$ (or simply $\Gamma$) is said to be \textit{strongly
irreducible} if there do not exist any finite union of proper
subspaces which is $\Gamma$-invariant.

The first main theorem of this paper is as follows:
\begin{theorem}\label{main}
Let $\Gamma$ be a sub-semigroup of $\m,$ $d>1,$ such that $\Gamma$
contains a proximal and expanding element and the $\Gamma$-action
on $\mathbb R^d$ is strongly irreducible. Then the semigroup
$\Gamma$ acting on $\T^d$ has the ID-property, that is, every
infinite $\Gamma$-invariant subset of $\T^d$ is dense.
\end{theorem}

If $d=1,$ one has $M_{inv}(1,\Z)=\Z^*\subset\R^*.$ As said above,
the conclusion of Theorem \ref{main} is valid in this case too, if
and only if $\Gamma$ is non-lacunary, i.e., not contained in a
cyclic subgroup of $\R^*.$ For $d>1,$ the condition in Theorem
\ref{main} imply that $\Gamma$ is not contained in a finite
extension of an abelian subgroup of $\m;$ in particular, here
$\Gamma$ is non-abelian, hence the situation of \cite{B} is
excluded from our setting.

The first step in order to get Theorem \ref{main} is to study
infinite $\Gamma$-invariant subsets $\Sigma$ of $\T^d$ such that 0
is a limit point of $\Sigma.$ Then we notice that the inverse
image in $\R^d$ of such an infinite $\Gamma$-invariant subset
contains some asymptotic set which consists of lines. Moreover,
there are some rays with good properties, that is which are not
contained in a subspace having a basis which consists of integer
vectors. This allows us to project them using canonical projection
$p:\R^d\to \R^d\slash\Z^d,$ $p(x)=x+\Z^d$ on $\mathbb T^d$ and
obtain the result in case when 0 is a limit point of the subset
$\Sigma$. Furthermore, using arguments close to \cite{B} and
\cite{F} and reduction to a finitely generated semigroup of
$\Gamma,$ we show that the opposite situation does not exist.

Let $L_\Gamma\subset\pr^{d-1}$ be the closure of the set of
directions corresponding to dominant eigenvectors of the proximal
elements in $\Gamma.$ We denote by $\tilde L_\Gamma$ the set of
corresponding nonzero vectors in $V=\R^d,$ by $\sigma$ the
symmetry $\sigma: v\mapsto -v$ in $V,$ and by $\tilde V$ the
factor space $\tilde V=V\slash\{\sigma,\text{Id}\}.$

The following is the basic tool in the proof of Theorem
\ref{main}.
\begin{theorem}\label{maintool}
Suppose $\Gamma$ is a sub-semigroup of $GL(d,\R),$ $d>1$ which is
strongly irreducible and contains a proximal and expanding
element. Let $\Sigma$ be a $\Gamma$-invariant subset of $\tilde
V\setminus\{0\}$ and suppose $0\in\overline{\Sigma}.$ Then
\begin{equation*}
\overline{\Sigma}\supset\tilde L_\Gamma\slash\{\sigma,\text{Id}\}.
\end{equation*}
\end{theorem}

To have a simple example in mind illustrating Theorem \ref{main}
and Theorem \ref{maintool}, consider the torus $\T^2.$ One of the
simplest example of a sub-semigroup of $SL(2,\Z)$ satisfying
properties of Theorem \ref{main} is the semigroup $\Gamma=\langle
a,b\rangle$ generated by two matrices $a=\bigl(\begin{smallmatrix}
2&1\\
1&1
\end{smallmatrix}\bigr)$
and $b=\bigl(\begin{smallmatrix}
3&2\\
1&1
\end{smallmatrix}\bigr)$
from $SL(2,\Z).$

>From Theorem \ref{main} we obtain that the $\Gamma$-orbits in
$\T^2$ are finite or dense. Furthermore we observe that, in the
context of Theorem \ref{maintool}, the dynamics of $\Gamma$ on
$\R^2$ is easy to visualize. The closure of the eigen-directions
in the positive quadrant $\R^2_+$ form a Cantor set and the
corresponding lines form an "attractor set" $\tilde L_\Gamma^1$
for the action of $\Gamma$ in $\R^2_+.$ There exist in $\R^2$
vectors whose orbit closures contain $0,$ for example dominant
eigenvectors of elements of $\Gamma^{-1}.$ The $\Gamma$-orbit for
such a vector tends to fill $\tilde L_\Gamma^1\cup -\tilde
L_\Gamma^1$ since the dynamics of its $\Gamma$-orbit consists of
attraction towards $0$ and expansion along the eigenvectors
sitting in $\tilde L_\Gamma^1\cup -\tilde L_\Gamma^1.$

For a general vector, for example a vector $v\in\R^2_+,$ there is
attraction towards $\tilde L_\Gamma^1$ and expansion along $\tilde
L_\Gamma^1$ and the $\Gamma$-orbit of $v$ is "thick at infinity"
due to the irrationality properties of eigenvalues of elements in
$\Gamma.$ In general the situation is similar, in particular the
projections of general $\Gamma$-orbits in $\T^d$ are large, hence
one can expect the closed $\Gamma$-orbits in $\T^d$ to be finite
unions of special manifolds, as conjectured in \cite{Ma}.

Let us consider now for $c>1,$ the factor space $\pr_c^{d-1}$ of
$V\setminus\{0\}$ by the subgroup of homotheties with ratio $\pm
c^k$ ($k\in\Z$). The action of $g\in G=GL(d,\R)$ on
$v\in\pr_c^{d-1}$ will be denoted $v\mapsto g.v.$ Let $\mu$ be a
probability measure on $\Gamma\subset GL(d,\R)$ whose support
generates $\Gamma.$ Then we can define an associated Markov
operator $P_\mu$ on $\pr_c^{d-1}$ by the formula:
\begin{equation*}
P_\mu(v,\cdot)=\int\delta_{g.v}d\mu(g).
\end{equation*}
The iterates $P_\mu^n$ of $P_\mu$ define a random walk on
$\pr_c^{d-1}.$

The following describes the asymptotic behavior of the iterates
$P_\mu^n;$ it is an essential tool in the proof of Theorem
\ref{maintool}, hence of Theorem \ref{main}.

\begin{theorem}\label{essentialtool}
Assume that $\Gamma\subset GL(d,\R)$ is a sub-semigroup which is
strongly irreducible and contains a proximal element. Then with
the above notations, the Markov operator $P_\mu$ on $\pr_c^{d-1}$
has a unique stationary measure $\rho,$ the support $S_\rho$ of
$\rho$ is the unique closed $\Gamma$-invariant minimal subset of
$\pr_c^{d-1}$ and for any $v\in\pr_c^{d-1}$ the sequence of
measures $P_\mu^n(v,\cdot)$ converges to $\rho.$ Moreover, the
trajectories of $P_\mu,$ starting from $v$ converge a.e. to
$S_\rho.$
\end{theorem}

Along the way, we get some new results and informations. For
example, we show a-priori that \textit{weak} ID-property (that is
the closures of the orbits $\Gamma x,$ $x\in\T^d$ are either
finite or equal $\T^d$ itself) and ID-property are equivalent, a
fact implicitly used in previous papers, but apparently unproved
in the literature.

We also clarify the relations between a fundamental cocycle
equation on $\Gamma\times\pr^{d-1}$ and an aperiodicity condition
for the dominant eigenvalues of proximal elements in $\Gamma$
which occur in \cite{K} and which has also a geometric
interpretation in terms of lengths of closed geodesics (see
\cite{Eb}).

Furthermore, the result in Theorem \ref{essentialtool} extends
results of \cite{GR} but is new in this generality.

Also the result of Theorem \ref{main} is not covered by \cite{GS}
since, in our setting, $\Gamma$ is allowed to be a sub-semigroup
of $\m,$ ($d>1$). Then we need to prove that $\Gamma$ can be
supposed to be finitely generated (see Proposition \ref{prop24}).

The structure of the paper is as follows. In Sect.
\ref{preliminaries} we set the notation and give all necessary
definitions. In particular, we define a dominant vector, a
proximal element and state our two hypothesis $(H_1)$ - strong
irreducibility and $(H_2)$ - proximality, under which we prove
Theorem \ref{essentialtool}. We introduce hypothesis $(H_0),$ i.e.
the unboundedness of $\Gamma$-orbits in $V\setminus\{0\}.$ Under
$(H_1)$ and $(H_2),$ this condition is equivalent (see Proposition
\ref{dodatkowa}) to the existence of a proximal and expanding
element in $\Gamma$ which allows to prove Theorem \ref{main} and
Theorem \ref{maintool}. It is clear, that this condition is
necessary for validity of ID-property.

We observe that conditions $(H_0),$ $(H_1)$ and $(H_2)$ are
analogous to those used in \cite[Theorem 2.5]{GLP} in order to get
a homogeneous behavior at infinity of the potential measure in
$\tilde V$ associated with $\mu,$ hence also of the
$\Gamma$-orbits at infinity in $\tilde V.$ (See also \cite{SGP}
for the case of affine actions.)

In Sect. \ref{equiv} we prove the equivalence of the weak
ID-property and ID-property (Proposition \ref{e}).

In Sect. ~\ref{asset} we study the $\Gamma$-actions on various
spaces, namely on the projective space $\pr(V),$ the compact
homogeneous space $\pr_c(V)$ and $V$ itself. We define the
asymptotic sets for $\Gamma$-actions and we study their
properties. We also clarify the role of aperiodicity hypotheses of
$\Gamma$ considered by Kesten in \cite{K} (see Corollary
\ref{compcoll}).

Sect. \ref{orbits} develops the random walks techniques which are
used in the proof of the main result of this section which is
Theorem \ref{tm}. This theorem together with the method presented
in \cite{F} allow us to prove in Sect. ~\ref{proof} Theorem
\ref{main}. Theorem \ref{tm} follows from a detailed study of
random walks on $V$ and various $\Gamma$-spaces, governed by a
measure $\mu$ sitting on $\Gamma$ and such that the convolution
iterations $\mu^{*k}$ fill $\Gamma.$ Some of these results are
well known but we have included the proofs in order to be
self-contained. Some others are new.

Finally, in Sect. ~\ref{proof} we give the proof of Theorem
\ref{main}.

\section{Proximality, irreducibility, expansivity}\label{preliminaries}
In what follows $\Gamma$ will denote a sub-semigroup of
$GL(d,\mathbb R).$ We consider the actions of $\Gamma$ on the
vector space $V=\mathbb R^d,$ on the associated projective space
$\mathbb P^{d-1}=\pr(V),$ and on $\tilde V=V\slash
\{\text{Id},\sigma\}=V\slash\{\pm\text{Id}\}.$ We denote by $\pi$
the projection of $V\setminus\{0\}$ on $\mathbb P^{d-1}=\pr(V)$
and we identify $\pr(V)$ with the unit sphere $\mathbb S^{d-1}$
divided by the symmetry $\sigma: x\mapsto -x.$ Also $K=SO(d,\R)$
will denote the special orthogonal group and $m$ the unique
$K$-invariant probability measure on $\pr(V).$

The action of the matrix $g$ on the vector $x\in V$ we denote by
$(g,x)\mapsto gx$, whereas for the action of $g$ on the projective
space $\pr(V)$ we write: $g.\pi(x)=\pi(gx).$

A matrix $\gamma\in GL(d,\mathbb R)$ is said to be
\textit{proximal} if it has an eigenvalue $\lambda_\gamma$ such
that $|\lambda_\gamma|>|\lambda|$ for all other eigenvalues
$\lambda$ of $\gamma.$ Thus $\lambda_\gamma\in\R.$ For such a
$\gamma$ an eigenvector $v_\gamma\in V$ corresponding to the
eigenvalue $\lambda_\gamma$ is called a \textit{dominant
eigenvector} or simply \textit{dominant vector} of $\gamma$. By
$\Delta_\Gamma$ we denote the set of all proximal elements in
$\Gamma.$ An element $\gamma\in GL(d,\R)$ is said to be
\textit{expanding} if it has an eigenvalue $\lambda$ such that
$|\lambda|>1.$

More generally, for $u\in\text{End}(V),$ we denote $|\lambda_u|$
the spectral radius of $u.$

If $\gamma\in\Delta_\Gamma$ then we define $\gamma^+\in\pr(V)$ as
a point corresponding to the line in $V$ generated by $v_\gamma.$
By $V^{<}_\gamma$ we denote the unique $\gamma$-invariant
hyperplane complementary to $V^{\text{max}}_\gamma=\R v_\gamma.$

We consider the following assumptions.
\begin{itemize}
\item [$(H_0)$] For every $v\in V\setminus\{0\},$ the orbit $\Gamma v$ is unbounded.
\item [$(H_1)$] The $\Gamma$-action is \textit{strongly
irreducible} (in short, $\Gamma$ is \textit{strongly
irreducible}), that is that there does not exist any finite union
of proper subspaces which is $\Gamma$-invariant.
\item [$(H_2)$] $\Gamma$ contains a matrix $\gamma$ which is proximal.
\end{itemize}
\begin{remark}\label{requiv}\emph{(i) The condition $(H_1)$ can be equivalently
formulated as follows. A sub-semigroup $\Gamma$ of $GL(V)$ acts
strongly irreducibly on $V$ if every finite index subgroup $H$ of
the group $\langle\Gamma,\Gamma^{-1}\rangle$ acts irreducibly on
$V,$ that is every $H$-invariant subspace of $V$ is either 0 or
$V.$}

\emph{(ii) If $\Gamma$ is a sub-semigroup of $SL(d,\R),$ then
conditions $(H_1)$ and $(H_2)$ imply $(H_0),$ since otherwise the
determinant of the proximal element $\gamma$ would be strictly
less than $1$ (see Proposition \ref{dodatkowa} b) below).}

\emph{(iii) The condition $(H_1)$ (resp. $(H_2)$) if valid for
$\Gamma,$ is also valid for $\Gamma^t,$ the transposed semigroup
acting on the dual space $V^*.$}

\emph{(iv) Conditions $(H_0),$ $(H_1)$ for $\Gamma$ imply
condition $(H_0)$ for $\Gamma^t.$ This will be used in the proof
of Theorem \ref{main} and can be seen as follows. Let $W\subset
V^*$ be the subspace of vectors with bounded $\Gamma^t$-orbits.
Then $W$ is $\Gamma^t$-invariant, hence (iii) implies $W=\{0\}$ or
$W=V^*.$ In case $W=V^*,$ $\Gamma^t$ is relatively compact in
$\text{End}(V^*),$ hence $\Gamma$ is relatively compact in
$\text{End}(V).$ This contradicts condition $(H_0)$ for $\Gamma.$}
\end{remark}

The concept of Zariski closure, defined below, will be freely used
in dealing with the above conditions (see \cite{OV}).

Let $\Gamma$ be a subset of $GL(d,\R).$ We recall that the
\textit{Zariski closure} $\text{Zc}(\Gamma)$ of $\Gamma$ is the
set of zeros of all real polynomials with variables in the
coefficients of $g\in GL(d,\R)$ and $(\det g)^{-1},$ which vanish
on $\Gamma.$

If $\Gamma$ is a sub-semigroup of $GL(d,\R)$ then
$\text{Zc}(\Gamma)$ is a group which contains $\Gamma,$ is closed
and has a finite number of connected components in the real
topology (see \cite{OV}). The connected component of the identity
in the Zariski topology is a subgroup of finite index which will
be denoted by $\text{Zc}_0(\Gamma).$

We have the following generalization of Lemma 2.8 in \cite{CG00}
to the case of semigroups.
\begin{lemma}\label{l2.8}
Let $\Gamma\subset GL(V)$ be a sub-semigroup. The $\Gamma$-action
satisfies condition $(H_1)$, if and only if there do not exist a
non-zero vector $v$ such that the orbit $\Gamma v$ is contained in
a finite union of proper vector subspaces of  $V.$
\end{lemma}
\begin{proof} Suppose $(H_1)$ to be valid and $v\in V$ be such
that $\Gamma v \subset\bigcup_{j=1}^nV_j$ where $V_j$ are proper
subspaces of $V.$ Let $W$ be a finite union of subspaces of $V$
such that $\Gamma v \subset W$ and $\mathcal{W}$ the set of such
$W.$ We observe that $\Gamma v \subset\bigcap_{W\in\mathcal{W}}W.$
Since a strictly decreasing family of elements of $\mathcal{W}$ is
finite  we get that $\bigcap_{W\in\mathcal{W}}W$ belongs also to
$\mathcal{W},$ in other words $W_0:=\bigcap_{W\in\mathcal{W}} W$
is the minimum element in $\mathcal{W}.$ We write
$W_0=\bigcup_{j=1}^mV_j$ and we are going to show that $W_0$ is
preserved by $\Gamma.$ Since $W\in\mathcal{W}$ is algebraic we
have: $\text{Zc}(\Gamma)v\subset W,$ in particular
$\langle\Gamma,\Gamma^{-1}\rangle v\subset W.$ It follows that,
for any $\gamma\in\Gamma:$
\begin{equation*}
\gamma W\supset\gamma\langle\Gamma,\Gamma^{-1}\rangle
v\supset\Gamma v.
\end{equation*}
Hence, $\gamma W\in\mathcal{W}.$ Since $W_0$ is the minimum
element of $\mathcal{W},$ we have $\gamma W_0\supset W_0,$ $\gamma
W_0=W_0;$ hence, $\Gamma W_0=W_0.$ Condition $(H_1)$ says that it
is impossible.

Conversely, suppose $V_j\,(1\leq j\leq n)$ is a family of proper
subspaces which is preserved by $\Gamma.$ Let $v\in V_1,$ then
$\Gamma v\in\bigcup_{i=1}^nV_i.$ From the hypothesis this is
impossible, hence condition $(H_1)$ is satisfied.
\end{proof}

Let $X$ be a compact metric space with distance function $\delta.$
We say that the action of a semigroup $\Gamma$ of continuous
transformations of $X$ is \textit{proximal} if, given $x,y\in X,$
there exists a sequence $\{\gamma_n\}\subset\Gamma$ such that
$\delta(\gamma_n.x,\gamma_n y)\to 0$ as $n\to\infty.$

Define the distance function $\delta$ on $\pr(V)$ as follows:
\begin{equation*}
\delta(\bar{u},\bar{v})=\|u\wedge v\|\slash \|u\|\|v\|,\;
\bar{u},\bar{v}\in\pr(V),
\end{equation*}
where $u$ and $v$ are the corresponding vectors in the vector
space $V.$

\begin{proposition}[Theorem 2.9 in \cite{GG}]\label{prox onpr}
Let $\Gamma$ be a sub-semigroup of the group $GL(V).$

Then we have the following equivalence:
\begin{itemize}
\item [(a)] $\Gamma$ satisfies $(H_1)$ and $(H_2).$
\item [(b)] $\Gamma$ acts proximally on $\pr(V)$ and is strongly
irreducible.
\end{itemize}
\end{proposition}
\begin{proof}
((a) $\Rightarrow$ (b)) We consider a proximal element
$\gamma\in\Gamma$ and denote
\begin{equation*}
\delta=\lim_n\|\gamma^{2n}\|^{-1}\gamma^{2n},\qquad
\mathfrak{z}=\text{Ker}\,\delta\subset\pr(V^*).
\end{equation*}
Then if $x,y\in\pr(V)$ do not belong to $\text{Ker}\,\delta,$ we
have $\lim_n\gamma^n.x=\gamma^+,$ $\lim_n\gamma^n.y=\gamma^+.$
Hence, $\lim_n\delta(\gamma^n.x,\gamma^n.y)=0.$

In general, if $x,y\in\pr(V)$ are given we can find $h\in\Gamma$
such that $h.x\not\in\text{Ker}\,\delta$ and
$h.y\not\in\text{Ker}\,\delta,$ otherwise, passing to the dual
space $V^*,$ transposing maps, and using the hyperplanes $x^\bot$
and $y^\bot$ of $V^*$ defined by $x$ and $y,$ one would have
\begin{equation*}
\forall\, h\in\Gamma:\;h^t.\mathfrak{z}\subset x^\bot\text{ or
}h^t.\mathfrak{z}\subset y^\bot.
\end{equation*}
But Remark \ref{requiv} iii) and Lemma \ref{l2.8} say that this is
impossible under condition $(H_2).$

((b) $\Rightarrow$ (a)) It follows from proximality of $\Gamma$ on
$\pr(V)$ (see \cite{F1}) that, given a finite subset
$E\subset\pr(V),$ there exist a sequence $\{g_n\}\subset\Gamma$
and $x\in\pr(V)$ such that
\begin{equation*}
\forall\, y\in E:\;\lim_ng_n.y=x.
\end{equation*}
We consider a finite system $E=\{x_1,x_2,\ldots,x_{2d-1}\}$ of
$2d-1$ points in $\pr(V)$ such that any $d-$subsystem consists of
independent points.

We consider the linear maps $u_n=\frac{g_n}{\|g_n\|}$ and using a
convergent subsequence, we can assume suppose than $u_n$ converges
towards $u\in\text{End}(V),$ $\|u\|=1.$ We show that $u$ has rank
one.

>From the definition of $E$ it follows that, at least $d$ points of
$E$ do not belong to $\text{Ker}\,u.$ We replace these points, as
well as $x,$ by the corresponding unit vectors in $V,$ say $\tilde
x_1,\ldots,\tilde x_d,\tilde x.$ Then we obtain
\begin{equation*}
u\tilde x_i=\lambda_i\tilde x,\;1\leq i\leq d,
\end{equation*}
where $\lambda_i\not=0;$ the points $\{\tilde x_i\}$ form a basis
of $V,$ hence the rank of $u$ is one, i.e.,
$\dim\text{Ker}\,u=d-1.$ We can moreover suppose that
$\text{Im}\,u\not\subset\text{Ker}\,u,$ since otherwise we could
replace $g_n$ by $gg_n,$ where $g\in\Gamma$ satisfies
$\text{Im}\,gu=g(\text{Im}\,u)\not\subset\text{Ker}\,u$ and
$\text{Ker}\,gu=\text{Ker}\,u.$ The existence of $g\in\Gamma$ such
that $g(\text{Im}\,u)\not\subset\text{Ker}\,u=\text{Ker}\,gu$
follows from Lemma \ref{l2.8}.

Under this condition, $u$ is proportional to the projection on
$\text{Im}\,u,$ along the hyperplane $\text{Ker}\,u.$ In
particular, $u$ has a unique non-zero eigenvalue $\lambda.$ Since
the sequence $\frac{g_n}{\|g_n\|}-u$ converges to zero, we obtain
that for $n$ large, $\frac{g_n}{\|g_n\|}$ has also a unique
dominant eigenvalue close to $\lambda.$ The same is true for
$g_n,$ hence $\Gamma$ satisfies $(H_2).$
\end{proof}
\begin{proposition}\label{dodatkowa}
Let $\Gamma$ be a sub-semigroup of $GL(V).$ Then we have the
following equivalence:
\begin{itemize}
\item [(a)] $\Gamma$ satisfies $(H_1),$ $(H_2)$ and the element $\gamma$ in the condition
$(H_2)$ satisfies $|\lambda_\gamma|>1.$
\item [(b)] $\Gamma$ is unbounded and satisfies $(H_1)$ and $(H_2).$
\item [(c)] $\Gamma$ satisfies $(H_0),$ $(H_1)$ and $(H_2)$.
\item [(d)] $\Gamma$ satisfies $(H_1),$ $(H_2)$ and there exists
$\gamma\in\Gamma$ such that $|\lambda_\gamma|>1.$
\end{itemize}
\end{proposition}
\begin{proof}
((d) $\Rightarrow$ (b)) Let $\gamma\in\Gamma$ be an expanding
element in $\Gamma,$ hence $|\lambda_\gamma|>1.$ Then:
$\|\gamma^n\|\geq|\lambda_\gamma|^n.$ Hence
$\lim_n\|\gamma^n\|=\infty,$ i.e., $\Gamma$ is unbounded.

((b) $\Rightarrow$ (a)) We will use the basic \cite[Theorem
4.1]{AMS}, which allows to construct new proximal maps and, which
implies the following. If $\Gamma\subset GL(V)$ satisfies $(H_1)$
and $(H_2)$ there exists $\varepsilon>0,$ $r>1,$  and a finite
subset $M\subset\Gamma,$ such that, for any $g\in GL(V),$ there
exists $a\in M$ such that $ag$ is proximal, the distance in
$\pr(V)$ of $(ag)^+$ to $V_{ag}^<$ is at least $\varepsilon,$ and:
\begin{equation*}
|\lambda_{ag}|\geq r\|(ag)|_{V^<_{ag}}\|.
\end{equation*}
Since $\Gamma$ is unbounded, there exists a sequence
$\{\gamma_n\}\subset\Gamma$ such that
\begin{equation*}
\lim_n\|\gamma_n\|=\infty.
\end{equation*}
Using a subsequence of $\gamma_n$ we can suppose that, for some
$a\in M,$ $a\gamma_n$ is proximal, $(a\gamma_n)^+,$
($V_{a\gamma_n}^<$ resp.) converges to $x\in\pr(V)$
($W_n=V^<_{a\gamma_n}$ converges to the hyperplane $W$ of
$\pr(V),$ resp.). We have $x\not\in W,$ since the distance of
$(a\gamma_n)^+$ to $W_n$ is at least $\varepsilon.$ We can suppose
also that $\frac{a\gamma_n}{\|a\gamma_n\|}$ converges to
$u\in\text{End}(V)$ with $\|u\|=1.$ Clearly, $V$ is the direct sum
of the hyperplane $W$ and of the line generated by $x.$
Furthermore,
$|\lambda_u|=\lim_n\frac{|\lambda_{a\gamma_n}|}{\|a\gamma_n\|},$
$\|u|_W\|=\lim_n\frac{1}{\|a\gamma_n\|}\|(a\gamma_n)|_{W_n}\|.$
Since $u\not=0$ preserves the above direct sum we have
$|\lambda_u|>0.$ Then the condition $|\lambda_{a\gamma_n}|\geq
r\|(a\gamma_n)|_{W_n}\|$ implies:
\begin{equation*}
|\lambda_u|\geq r\|u|_W\|,\;|\lambda_u|>\|u|_W\|.
\end{equation*}
In particular, $u$ has a dominant eigenvalue which is simple.
Since, $\|u-\frac{a\gamma_n}{\|a\gamma_n\|}\|$ converges to zero,
we have for $n$ large:
\begin{equation*}
|\lambda_{a\gamma_n}|\geq\|a\gamma_n\|\frac{|\lambda_u|}{2}.
\end{equation*}
Moreover, the condition $\lim_n\|\gamma_n\|=\infty$ implies
\begin{equation*}
\lim_n\|a\gamma_n\|\geq\lim_n\|a^{-1}\|^{-1}\|\gamma_n\|=\infty.
\end{equation*}
In particular, for $n$ large, $|\lambda_{a\gamma_n}|>1,$ hence
$a\gamma_n$ is proximal and expanding, i.e., (a) is valid.

((b) $\Rightarrow$ (c)) We consider the subspace $W\subset V$ of
vectors in $V$ having a bounded $\Gamma$-orbit. Clearly, this
subspace is $\Gamma$-invariant. Then condition $(H_1)$ implies
$W=V$ or $W=\{0\}.$ In the second case $(H_0)$ has been proved.
The first case do not occur since it contradicts the hypothesis
that $\Gamma$ is unbounded.

((c) $\Rightarrow$ (b)) and ((a) $\Rightarrow$ (d)) are trivial.
\end{proof}

The following is a useful characterization of strong
irreducibility in terms of Zariski closure.
\begin{proposition}\label{prop23}
Let $\Gamma$ be a sub-semigroup of $GL(V).$ Then $\Gamma$
satisfies $(H_1)$ if and only if $\text{Zc}_0(\Gamma)$ acts
irreducibly on $V.$
\end{proposition}
\begin{proof}
Assume that $\Gamma$ satisfies $(H_1)$ and let $W\subset V$ be a
non-zero and $\text{Zc}_0(\Gamma)$-invariant subspace. For some
finite set $F\subset\Gamma$ we have
$\Gamma\subset\text{Zc}(\Gamma)=\bigcup_{\gamma\in
F}\gamma\text{Zc}_0(\Gamma),$ hence $\Gamma W=\bigcup_{\gamma\in
F}\gamma W.$ Since $\Gamma$ satisfies $(H_1)$ we get $W=V,$ hence
$\text{Zc}_0(\Gamma)$ acts irreducibly on $V.$

Assume that $\text{Zc}_0(\Gamma)$ acts irreducibly on $V$ and let
$W$ be a non-zero subspace of $V,$ $F$ a finite subset of $\Gamma$
such that $\Gamma W=\bigcup_{\gamma\in F}\gamma W.$ Since $\Gamma
W$ is an algebraic manifold, $\text{Zc}(\Gamma)$ leaves $\Gamma W$
invariant, hence permutes the subspaces $\gamma W$
($\gamma\in\Gamma$). Since $\text{Zc}_0(\Gamma)$ is connected, we
have for any $\gamma\in F:$
\begin{equation*}
\text{Zc}_0(\Gamma)\gamma W=\gamma W.
\end{equation*}
>From the irreducibility of the action of $\text{Zc}_0(\Gamma)$ on
$V,$ we get $W=V.$
\end{proof}
The following will be essential in the proof of Theorem \ref{main}
\begin{proposition}\label{prop24}
Assume that the semigroup $\Gamma\subset GL(V)$ satisfies $(H_1)$
and $(H_2).$ Then $\Gamma$ contains a finitely generated
sub-semigroup which satisfies $(H_1)$ and $(H_2).$
\end{proposition}
The proof of the above proposition depends on the following
\begin{lemma}\label{lem25}
Assume that $\Gamma$ satisfies $(H_1),$ $(H_2).$ Denote by $D$
(resp. $C$) the commutator subgroup (resp. connected center) of
$\text{Zc}_0(\Gamma).$ Then $\text{Zc}_0(\Gamma)$ is the almost
direct product of $D$ and $C.$ Furthermore, $D$ is semisimple
without compact factors and $C$ consists of homotheties.
\end{lemma}
\begin{proof}
Since $\Gamma$ acts irreducibly on $V,$ $\text{Zc}(\Gamma)$ is a
$\R$-reductive group (see \cite{OV}), hence $D$ is semisimple and
$\text{Zc}_0(\Gamma)$ is the almost direct product of $C$ and $D.$
We can write $D$ as as the almost direct product $D=D_1D_2,$ where
$D_1$ is compact and $D_2$ is semisimple without compact factor.
Since $\Gamma$ contains a proximal element and
$\text{Zc}(\Gamma)\slash\text{Zc}_0(\Gamma)$ is finite,
$\text{Zc}_0(\Gamma)$ contains also a proximal element. We denote
by $\gamma$ this element and write $\gamma v_\gamma=\lambda_\gamma
v_\gamma,$ $\gamma=cd_1d_2$ with $c\in C,$ $d_1\in D_1,$ $d_2\in
D_2.$ Since $d_1$ and $\gamma$ commute, and the direction of
$v_\gamma$ is uniquely defined by $\gamma,$ $d_1v_\gamma$ is
proportional to $v_\gamma.$ Since $D_1$ is compact we have
$d_1v_\gamma=\pm v_\gamma,$ hence $cd_2$ is also proximal with
dominant eigenvector $v_\gamma.$ Since $C$ commutes with $cd_2,$
and $v_\gamma$ is $cd_2$-dominant, there exists a $\R$-character
$\chi$ of $C$ such that for every $g\in C:$
$gv_\gamma=\chi(g)v_\gamma.$ Since the subspace
$W=\{v:cv=\chi(c)v,\forall c\in C\}$ is $\Gamma$-invariant,
contains $v_\gamma$ and the action of $\Gamma$ is irreducible, we
have that for every $v\in V$ and for every $c\in C,$
$cv=\chi(c)v.$ Thus, $C$ consists of homotheties, $D_1D_2$ acts
also irreducibly on $V$ and $v_\gamma$ is $d_2$-dominant. Since
$D_1$ commutes with $d_2$ we get, as above, that $D_1$ preserves
the direction of $v_\gamma.$ Since $D_1$ is compact and connected,
we obtain that $v_\gamma$ is $D_1$-invariant. Since $D_1$ commutes
with $CD_2,$ the subspace of $D_1$-invariant vectors is preserved
by the action of $CD_1D_2.$ From the irreducibility of
$\text{Zc}_0(\Gamma),$ we get that $D_1=\text{Id},$ hence
$\text{Zc}_0(\Gamma)=CD_2.$
\end{proof}
\begin{proof}[Proof of Proposition \ref{prop24}]
We consider the semigroup $\Gamma(S)$ generated by the finite set
$S\subset\Gamma.$ Clearly, if $S\subset S^\prime,$ then
$\Gamma(S)\subset\Gamma(S^\prime).$ We take a totally ordered
family $S_i$ ($i\in I$) such that $\Gamma=\bigcup_{i\in
I}\Gamma(S_i);$ we denote by $G_0^i$ the connected component of
the identity in $\text{Zc}(\Gamma(S_i)).$ Then, since
$G_0^i\subset G_0^j$ if $S_i\subset S_j,$ we get that for some
$\iota\in I:$
\begin{equation*}
\text{$H_0:=G_0^\iota=\bigcup_{i\in I}G_0^i=G_0^k$ if $S_k\supset
S_\iota.$}
\end{equation*}
We can suppose that for any $i\in I,$ $G_0^i=G_0^\iota.$ It
follows that $H_0$ is normal in $\text{Zc}(\Gamma(S_i))$ for any
$i\in I,$ hence $H_0$ is normal in $\text{Zc}(\Gamma).$ In
particular, $H_0$ is contained in $\text{Zc}_0(\Gamma).$ We
observe that $H_0$ has finite index in $\text{Zc}(\Gamma(S_i)),$
hence $L=\text{Zc}_0(\Gamma)\slash H_0$ is an algebraic group
which is the Zariski closure of the union of the finite subgroups
$\Phi_i$ corresponding to $\text{Zc}_0(\Gamma(S_i)).$ In view of
Lemma \ref{lem25} we know that the algebraic group $L$ has the
same structure as $\text{Zc}_0(\Gamma),$ in particular is
reductive. We write it as the almost direct product of its
connected center $C^\prime\subset\R^*$ and its commutator subgroup
$D^\prime.$ Passing to the factor group $L\slash D^\prime,$ using
the finite subgroups $\Phi_i,$ we get $C^\prime=\{\text{Id}\},$
$L=D^\prime.$ We consider a faithful, irreducible representation
of the adjoint group of $L$ into a real vector space $V^\prime.$
Then each finite subgroup $\Phi_i$ leaves a positive definite
quadratic form $q_i$ invariant. We can suppose that the forms
$q_i$ are normalized and we denote by $q$ a cluster value of the
$(q_i)_{i\in I}.$ Then $q$ is invariant under the action of
topological closure $\Phi$ of $\bigcup_{i\in I}\Phi_i.$ Since
$\text{Zc}(\bigcup_{i\in I}\Phi_i)=L=\text{Zc}(\Phi),$ we get that
$\Phi$ acts irreducibly on $V^\prime,$ as $L.$ Since the kernel of
$q$ is $\Phi$-invariant, we get that this kernel is trivial, hence
$q$ is positive definite. It follows that $\Phi$ is compact. Since
$\text{Zc}(\Phi)=L,$ we conclude that $L=\Phi$ is compact, hence
from Lemma \ref{lem25}, $L=\{\text{Id}\}.$ It follows:
$H_0=\text{Zc}_0(\Gamma)=\text{Zc}_0(\Gamma(S_\iota)).$ We can
suppose that $\Gamma(S_\iota)$ contains a proximal element from
$\Gamma\cap\text{Zc}_0(\Gamma).$ Then $\Gamma(S_\iota)$ is
finitely generated, and satisfies $(H_2).$ The condition $(H_1)$
is also satisfied by $\Gamma(S_\iota),$ since
$\text{Zc}_0(\Gamma(S_\iota))=\text{Zc}_0(\Gamma)$ acts
irreducibly on $V.$
\end{proof}
\begin{remark}\emph{
We will see in Lemma \ref{l1} below that condition $(H_0)$ remains
also valid while passing to a convenient finitely generated
sub-semigroup. However, this property cannot be achieved with
condition $(H_1)$ alone. A simple counterexample is the following:
suppose $\Gamma$ is the semigroup of rational rotations of the
Euclidean plane, centered at the origin. Then any finitely
generated sub-semigroup $\Gamma^\prime$ preserves a regular
polygon inscribed in the unit circle. Hence, condition $(H_1)$ is
not satisfied by $\Gamma^\prime.$}

\emph{This explains why we consider $(H_1)$ and $(H_2)$
simultaneously in Proposition \ref{prop24}.}
\end{remark}
\section{Equivalence of the weak ID-property and
ID-property}\label{equiv} Let us recall the definitions of the
weak ID-property and ID-property  once again in the context of
sub-semigroups of $\m$ acting in the usual way  on $d$-dimensional
tori. We say that a sub-semigroup $\Gamma$ of $\m$ has
\textit{ID-property} if every infinite $\Gamma$-invariant subset
of $\T^d$ is dense in $\T^d.$ This is of course equivalent to the
fact that every infinite closed $\Gamma$-invariant subset of
$\T^d$ is $\T^d$ itself.

Let us now recall the definition of the weak ID-property which is
defined in terms of orbits. We say that a sub-semigroup $\Gamma$
of $\m$ acting on $d$-dimensional torus has the \textit{weak
ID-property} if for every $x\in\T^d,$ the closure of the orbit
$\overline{\Gamma x}$ is either finite or the whole $\T^d.$

Here it is convenient to use condition ($H_0$) which is weaker
than the hypothesis in Theorem \ref{main}.
\begin{proposition}\label{e} Let $\Gamma$ be a sub-semigroup of
$\m$ acting on $\T^d$ and satisfying assumption $(H_0).$ Then
$\Gamma$ has weak ID-property if and only if $\Gamma$ has
ID-property.
\end{proposition}
Before we go to the proof of the above equivalence we need the
following three lemmas.
\begin{lemma}\label{l} Suppose $S$ is a finite subset of $GL(d,\R),$
which generates a semigroup $\Gamma$ which satisfies $(H_0).$
Denote
\begin{equation*}
C=\sup\{\|s\|:s\in S\}.
\end{equation*}
Then for every $x\in V,$ $\|x\|\leq 1$ there exists an element
$g\in\Gamma$ such that
\begin{equation*}
1<\|gx\|\leq C.
\end{equation*}
\end{lemma}
\begin{proof}
Note that $C>1,$ since $\Gamma$ is unbounded. We consider a
sequence $s_k\in S$ such that the sequence $s_n\ldots s_1x,$
$\|x\|\leq 1,$ is unbounded and denote
\begin{equation*}
k=\sup\{n\in\N:\,\|s_n\ldots s_1x\|\leq 1\}.
\end{equation*}
Then we have
\begin{equation*}
 \|s_{k+1}\|\leq C,\;\|s_k\ldots s_1x\|\leq 1,\;\|s_{k+1}\ldots
 s_1x\|>1.
\end{equation*}
It follows
\begin{equation*}
1<\|s_{k+1}s_k\ldots s_1x\|\leq C\|s_k\ldots s_1x\|\leq C.
\end{equation*}

Then the conclusion follows with $g=s_{k+1}\ldots s_1.$
\end{proof}
\begin{lemma}\label{l1}
Let $\Gamma$ be a sub-semigroup of $GL(d,\R)$ which satisfies
condition $(H_0).$ Then $\Gamma$ contains a finitely generated
sub-semigroup which satisfies $(H_0).$
\end{lemma}
\begin{proof}
For any finite subset $S\subset\Gamma,$ we denote by $\Gamma(S)$
the semigroup generated by $S,$ by $V(S)$ the subspace of vectors
$v\in V$ such that $\Gamma(S)v$ is bounded. We observe that the
inclusion $S\subset S^\prime$ implies $V(S^\prime)\subset V(S)$ We
consider a totally ordered family of finite subsets
$S_\iota,$($\iota\in I$) such that $\Gamma=\bigcup_{\iota\in I}
\Gamma(S_\iota).$ Then $W=\bigcap_{\iota\in I}V(S_\iota)$ is of
the form $V(S_j)$ for some $j\in I$ and we have
$V(S_\iota)=V(S_j)$ if $S_\iota\supset S_j.$ It follows that $W$
is $\Gamma$-invariant. Furthermore, for any $v\in W,$ $\iota\in
I,$ $\Gamma(S_\iota)v$ is bounded.

We show that, if $W\not=\{0\},$ then $\Gamma v$ is bounded, for
some $v\in W\setminus\{0\}.$ Hence using condition $(H_0):$
$W=\{0\}.$ This implies that condition ($H_0)$ is satisfied by
$\Gamma(S_j).$

We consider the complexified vector space $W^\C\subset V^\C,$ a
$\Gamma$-irreducible subspace $U\subset W^\C,$ and the action of
$\Gamma(S_\iota)$ on $U.$ Since every $\Gamma(S_\iota)$-orbit in
$W$ is bounded for any $\gamma\in\Gamma(S_\iota)$ we have
$|\lambda_\gamma|\leq 1,$ hence: for any $\gamma\in\Gamma,$
$|\tr\gamma|_U|\leq\dim U.$

Since the action of $\Gamma$ on $U$ is irreducible, Burnside's
theorem implies that the algebra $\text{End}(U)$ is generated by
$\Gamma$, i.e. there exist $\gamma_1,\ldots,\gamma_r$ in $\Gamma$
such that the linear forms $f_1,\ldots, f_r$ on $\text{End}(U)$
defined by
\begin{equation*}
f_k(w)=\text{Tr}\,(\gamma_kw)\; (1\leq k\leq r)
\end{equation*}
form a basis of $(\text{End}(U))^*.$ It follows that for every
$\gamma\in\Gamma,$ $|f_k(\gamma)|\leq\dim U.$ Since the $\{f_k\}$
form a basis of $(\text{End}(U))^*,$ we obtain that $\Gamma|_U$ is
bounded. Then any $\Gamma$-orbit in $U$ is bounded. Then the same
is true for the conjugated space $\bar{U}\subset V^\C,$ and for
the sum $U+\bar{U}\subset V^\C.$ In particular, any
$v\in(U+\bar{U})\cap V$ has a bounded $\Gamma$-orbit. Hence from
condition $(H_0),$ $(U+\bar{U})\cap V=\{0\},$ $U=\{0\},$
$W=\{0\}.$
\end{proof}

Let $B_\varepsilon\subset \R^d$ denote the ball with the radius
$\varepsilon$ and the center 0. If we consider $\T^d$ and
$\varepsilon<1\slash 2,$ then we denote also by $B_\varepsilon$
the homeomorphic image of the ball $B_\varepsilon\subset \R^d$ by
the canonical quotient map $p:\R^d\to\T^d=\R^d\slash\Z^d.$
\begin{lemma}\label{ll34}
Let $\Gamma$ be a sub-semigroup of $\m$ which satisfies $(H_0).$
Then, there exists $\varepsilon=\varepsilon_\Gamma>0$ such that
for every $0\not=x\in\T^d:$
\begin{equation*}
\Gamma x\cap\T^d\setminus B_\varepsilon\not=\{0\}.
\end{equation*}
\end{lemma}
\begin{proof}
>From Lemmas \ref{l1} and \ref{l} above we can find $C>1$ such that
for any $x\in B_\varepsilon\subset\R^d,$ $\varepsilon<1\slash 2,$
there exists $g\in\Gamma$ such that
\begin{equation*}
\varepsilon\leq\|gx\|\leq C\varepsilon.
\end{equation*}
If $\varepsilon_\Gamma=1\slash 2C<1\slash 2$ we get that, for
every $x\in B_\varepsilon\subset\T^d,$ we have $\Gamma
x\not\subset B_\varepsilon.$ If, for some $y\not\in B_\varepsilon$
we had $\Gamma y\subset B_\varepsilon,$ then, for some
$\gamma\in\Gamma,$ $x=\gamma y\in B_\varepsilon;$ hence, from the
above observation $\Gamma x\not\subset B_\varepsilon;$ in
particular, since $\Gamma x\subset\Gamma y,$ we have $\Gamma
y\not\subset B_\varepsilon$ and this is a contradiction.
\end{proof}
Now we are ready to give the proof of Proposition \ref{e}.
\begin{proof}[Proof of Proposition \ref{e}]
It is obvious that ID-property implies weak ID-property. Therefore
we have to prove the converse.

If for some $x\in\Sigma,$ $\Gamma x$ is infinite then the
hypothesis implies $\overline{\Gamma x}=\T^d=\Sigma.$ Hence we can
suppose that $\Sigma$ is infinite,
$\Sigma=\bigcup_{x\in\Sigma}\Gamma x$ and each $\Gamma x,$
$x\in\Sigma$ is finite. It follows that $\Sigma\subset p(\Q^d),$
hence $\Sigma$ is countable.

Now consider the sequence of derived sets,
\begin{equation}\label{seq}
\Sigma^0=\Sigma\supset\Sigma^1\supset\ldots\supset\Sigma^n\supset\ldots,
\end{equation}
that is $\Sigma^{n+1}$ is the set of limit points of $\Sigma^n.$
Actually, the sequence \eqref{seq} terminates, i.e. there is an
index $n$ such that $\Sigma^n=\emptyset.$ If not we consider
$\Sigma^\infty:=\bigcap_{n=0}^\infty\Sigma^n.$ Clearly,
$\Sigma^\infty$ is a closed and countable set such that the set
$(\Sigma^\infty)^{\text{ac}}$ of limit points of $\Sigma^\infty$
is equal to $\Sigma^\infty.$
 This means that we
have so called \textit{perfect set} $\Sigma^\infty$ which is
non-void and countable. The Baire theorem says that, since every
point is closed and has empty interior in $\Sigma^\infty,$ the set
$\Sigma^\infty$ has the same property, i.e. has empty interior in
$\Sigma^\infty,$ which is impossible. Therefore, there is $n\in\N$
such that
\begin{equation*}
\Sigma^0=\Sigma\supset\Sigma^1\supset\ldots\supset\Sigma^n=\emptyset.
\end{equation*}
Without loss of generality we may assume $n=2.$ It follows that
$\Sigma^1$ is finite. In fact, otherwise $\Sigma^2$ would not be
an empty set. Let $\{x_1,\ldots x_n\}=\Sigma^1\subset p(\Q^d)$ be
the set of limit points of $\Sigma$ and let $q$ be a common
denominator of $x_i,$ $1\leq i\leq n.$ Then $0\in q\Sigma^1$ is
the unique limit point of $q\Sigma.$ Let us consider a ball
$B_\varepsilon$ around $0$ with $\varepsilon<\varepsilon_\Gamma$
given by Lemma \ref{ll34}. Then the points of $q\Sigma$ outside
$B_\varepsilon$ have no limit point, hence form a finite set $F.$
Now we can consider the $\Gamma$-orbits of these points, i.e.
$\Gamma x,$ $x\in F.$ They form a $\Gamma$-invariant finite set
$F^\prime=\bigcup_{x\in F}\Gamma x$ that we can exclude from
$q\Sigma$ without changing its properties. Therefore, now we have
the new set $\Sigma^\prime=q\Sigma\setminus F^\prime$ which is
closed, infinite, $\Gamma$-invariant and fully included in
$B_\varepsilon,$ and this contradicts Lemma \ref{ll34}.
\end{proof}

\section{Asymptotic sets for $\Gamma$-actions and a cohomological equation} \label{asset}

As in Sect. \ref{preliminaries}, $\Gamma$ is a sub-semigroup of
$GL(V)=GL(d,\mathbb R)$ and we consider its action on $V$ and
$\pr(V)=\pi(V\setminus\{0\}).$

We define the \textit{asymptotic sets}:
\begin{equation*}
\begin{split}
L_\Gamma=&\overline{\{\pi(v_0):v_0\text{ is a
dominant vector for $\Gamma$}\}}\\
\tilde L_{\Gamma}=&\{v\not =0:\pi(v)\in
L_\Gamma\}=\pi^{-1}(L_\Gamma).
\end{split}
\end{equation*}
We start with the following proposition
which is a semigroup version of a result of Guivarc'h and Conze
(cf. \cite{CG00} Proposition 2.2).
\begin{proposition}\label{rgs}
Let $\Sigma$ be a $\Gamma$-invariant subset of $V\setminus\{0\}$
such that $0\in\overline{\Sigma}.$ Then, under assumptions
$(H_0),$ $(H_1)$ and $(H_2),$ for any proximal and expanding
element $\gamma\in\Gamma$ there exists a $\gamma$-dominant vector
$u_0$ such that
\begin{equation}\label{gz}
\gamma^{\Z}u_0:=\{\gamma^ku_0:k\in\Z\}\subset\overline{\Sigma.}
\end{equation}
\end{proposition}
\begin{proof}
Let $V=V_\gamma^{\text{max}}\oplus V_\gamma^<$ be the
decomposition of the space $V$ relative to a proximal and
expanding element $\gamma\in\Gamma.$

Let $x_i\in\Sigma$ and $x_i\to 0$ as $i\to\infty.$ Then there
exists a sequence $\{\alpha_i\}$ of reals and $w\in V$ such that
$\alpha_ix_i\to w$ as $i\to\infty.$ We will show that without loss
of generality we can assume that $w\not\in V_\gamma^<.$ In fact,
since $\Gamma$ acts strongly irreducibly on $V,$ by Lemma
\ref{l2.8}, one can find an element $h\in\Gamma$ such that
$hw\not\in V_\gamma^<,$ i.e. $w\not\in h^{-1}V_\gamma^<.$ Define
\begin{equation*}
\Gamma_1=h^{-1}\Gamma h,\, \gamma_1=h^{-1}\gamma h\in\Gamma_1,\,
\Sigma_1=h^{-1}\Sigma.
\end{equation*}
Then $\gamma_1$ is proximal in $\Gamma_1,$ $w\not\in
h^{-1}V_\gamma^<=V_{\gamma_1}^<,$ and $\Sigma_1$ is a
$\Gamma_1$-invariant subset that contains 0 as a limit point.
Assume that we have found a non-zero vector $u_0\in
V_{\gamma_1}^{\text{max}},$ such that
$\gamma_1^{\Z}u_0\subset\overline{\Sigma}_1.$ Then
$h^{-1}\gamma^{\Z}hu_0\subset h^{-1}\overline{\Sigma},$ i.e.,
$\gamma^{\Z}hu_0\subset\overline{\Sigma}.$ But $hu_0\in
V_\gamma^{\text{max}}$ and we are done.

Thus from the very beginning we can assume that $w\not\in
V_\gamma^<.$

Let $e_1,\ldots,e_n$ be a basis of $V$ such that $e_1\in
V^{\text{max}}_\gamma,$ $\|e_1\|=1$ and $e_2,\ldots,e_n\in
V_\gamma^<.$ Let $\phi_j:V\to\R$ be linear forms  such that
\begin{equation*}
x=\sum_{j=1}^n\phi_j(x)e_j,\; x\in V.
\end{equation*}
Let $\Phi:V=V_\gamma^{\text{max}}\oplus V^<_\gamma\to
V^{\text{max}}_\gamma$ be the projection along $V^<_\gamma,$ i.e.,
$\Phi(x)=\phi_1(x)e_1.$ Since $w\not\in V_\gamma^<,$ it follows
without loss of generality that $\Phi(x_i)\not=0$ (since
$\alpha_ix_i\to w\not\in V_\gamma^<$). Since
$|\lambda_\gamma|=\lambda>1,$ there exists a sequence of integers
$\{p_i\}$ such that $p_i\to\infty$ and
\begin{equation*}
1\leq\lambda^{p_i}|\phi_1(x_i)|\leq\lambda.
\end{equation*}
Now passing to a subsequence if necessary one can find a
$\gamma$-dominant vector $u_0\in V^{\text{max}}_\gamma$ such that
$\gamma^{p_i}(\Phi(x_i))\to u_0$ as $i\to\infty.$ We will prove
that $\gamma^{p_i}(x_i)\to u_0$ as $i\to\infty.$

Clearly, it is enough to show that
\begin{equation}\label{warunek}
\phi_j(x_i)\gamma^{p_i}(e_j)\to 0\text{ as $i\to \infty,$ for any
$j>1.$}
\end{equation}
In fact, $\gamma^{p_i}(x_i)\to u_0$ if
$\gamma^{p_i}(\sum_{j=1}^n\phi_j(x_i)e_j)=\sum_{j=1}^n\phi_j(x_i)\gamma^{p_i}(e_j)\to
u_0$ and so we are led to \eqref{warunek}.

Therefore, we are going to prove \eqref{warunek}. One has
\begin{equation*}
\phi_j(x_i)\gamma^{p_i}(e_j)=\frac{\phi_j(x_i)}{|\phi_1(x_i)|}\frac{\gamma^{p_i}(e_j)}{\lambda^{p_i}}
\|\gamma^{p_i}(\Phi(x_i))\|,
\end{equation*}
where the first fraction tends to $\phi_j(w)\slash|\phi_1(w)|,$
the second one tends to zero, and the third term tends to
$\|u_0\|$ as $i\to\infty.$

Since $x_i\in\Sigma$ and $\gamma^{p_i}\in\Gamma$ for every
$i\in\N$ it follows that $u_0\in\overline{\Sigma}.$ Since
$\lim_i\gamma^{p_i}(x_i)=u_0,$ we get that
$\gamma^{-m}(u_0)=\lim_i\gamma^{p_i-m}(x_i).$ We see that for
almost all $i,$ $\gamma^{p_i-m}(x_i)\in\Sigma,$ thus
$\gamma^{-m}u_0\in\overline{\Sigma}.$
\end{proof}
\begin{remark}\label{r0}
Notice that the condition \eqref{gz} implies that
$0\in\overline{\Sigma}.$ In fact, simply take a sequence $\Z\ni
k_n$ such that $k_n\to -\infty.$
\end{remark}
\begin{proposition}\label{4.1}
Under conditions $(H_1)$ and $(H_2),$ the set $L_\Gamma$ is the
unique minimal $\Gamma$-invariant closed subset of $\pr(V).$
\end{proposition}
\begin{proof} We first show
that $L_\Gamma$ is $\Gamma$-invariant. Consider
$g\in\Delta_\Gamma,$ and $u=\lim_n\frac{g^{2n}}{\|g^{2n}\|},$
where the limit exists in $\text{End}(V)$ by proximality of $g.$
Consider a decomposition $V=V_g^<\oplus V_g^{\text{max}}.$ Then
$u$ is a multiple of the projection of $V$ onto
$V_g^{\text{max}},$ along $V_g^<.$

On the other hand, we consider $\gamma\in\Gamma$ and want to show
that $\gamma.g^+\in L_\Gamma.$ We observe that for any
$\delta\in\Gamma,$ we have
\begin{equation*}
\lim_n\left\|\gamma\frac{g^{2n}}{\|g^{2n}\|}\delta-\gamma
u\delta\right\|=0.
\end{equation*}
We have $\text{Im}(\gamma u\delta)=\gamma(\text{Im}\,u)$ and
$\text{Ker}(\gamma u\delta)=\delta^{-1}(\text{Ker}\, u).$ We note
that $\gamma u\delta$ has rank one, like $u,$ hence $\gamma
u\delta$ will be a multiple of a one-dimensional projection if
$\gamma(\text{Im}\, u)\not\subset\delta^{-1}(\text{Ker}\, u),$
i.e. $\delta(\text{Im}(\gamma u))\not\subset\text{Ker}\, u.$

Since $\Gamma$ is strongly irreducible, Lemma \ref{l2.8} shows
that such a $\delta\in\Gamma$ exists. Then, having $\delta$ chosen
as above, we know from perturbation theory, as in the proof of
Proposition \ref{prox onpr}, that for $n$ large
$\gamma\frac{g^{2n}}{\|g^{2n}\|}\delta$ has a simple dominant
eigenvalue and the corresponding eigenvector is close to
$\gamma(v_g).$ In other words:
\begin{equation*}
\gamma g^{2m}\delta\in\Delta_\Gamma,\, \, \gamma.g^+=\lim_m(\gamma
g^{2m}\delta)^+.
\end{equation*}
Hence $\gamma L_\Gamma\subset L_\Gamma.$

Now let $\Lambda$ be a closed $\Gamma$-invariant subset of
$\pr(V)$ and let us show that $\Lambda\supset L_\Gamma.$ Since
$\Gamma$ is strongly irreducible, $\Lambda$ is not contained in a
proper subspace. In particular $\Lambda\not\subset V_\gamma^<,$
hence there exists $x\in\Lambda$ with $x\not\in V_\gamma^<.$ Then
$\gamma^+=\lim_n\gamma^n.x\in\Lambda.$ Since $\Lambda$ is closed,
we have $L_\Gamma\subset\Lambda.$ This shows that $L_\Gamma$ is
minimal and is the unique minimal subset of $\pr(V).$
\end{proof}
\begin{proposition}[Proposition 2.2 in \cite{CG01}]\label{2.2}
Let $\Gamma$ be  a sub-semigroup of the group $GL(d,\R),$ $d>1,$
satisfying hypothesis $(H_1)$ and $(H_2).$ Let $S$ be a subset of
$\Gamma$ generating it. If $\varphi$ is a non-zero, continuous
function on $\pr^{d-1},$ $t$ is real and $\theta\in[0,2\pi),$ then
the following equation:
\begin{equation*}
\forall\,\gamma\in S,\,\forall\, x\in L_\Gamma:\,
\varphi(\gamma.x)\|\gamma x\|^{it}=e^{i\theta}\varphi(x)
\end{equation*}
has no solution, unless $\theta=0,$ $t=0,$
$\varphi(x)\equiv\text{const}$ on $L_\Gamma.$
\end{proposition}
\begin{proof}
Clearly we can suppose that $|\det\gamma|=1$ for
$\gamma\in\Gamma.$ Consider the function $\psi(v)$ defined on
$\tilde L_\Gamma$ by formula $\psi(v)=\varphi(\pi(v))\|v\|^{-it}.$
Then the relation for $\varphi$ is
\begin{equation*}
\forall\, \gamma\in S,\, \psi(\gamma v)=e^{i\theta}\psi(v).
\end{equation*}
Suppose that $t\not=0$ and put $\log\rho=2\pi\slash |t|.$ Then
additionally we have $\psi(\pm\rho^kv)=\psi(v)$ and the condition
$\psi(\lambda v)=\psi(v)$ for some $v\in\tilde L_\Gamma$ and some
$\lambda\in\R^+_*$ implies that $\lambda=\pm\rho^k,$ where
$k\in\Z.$

Let $c$ be any of the values of $\psi$ and put
$L_c=\psi^{-1}(\{c\})\subset\tilde L_\Gamma.$ Then, since $\psi$
is continuous, $L_c$ is a nonempty closed subset of $\tilde
L_\Gamma,$ which satisfies
\begin{equation*}
\forall\, \gamma\in S,\, \gamma(L_c)\subset L_{ce^{i\theta}}.
\end{equation*}
We also have for every $\lambda\in H_\rho$, which is the group of
homotheties of the form $\pm\rho^k,$ $k\in\Z,$
\begin{equation*}
\lambda L_c=L_c.
\end{equation*}
If now $u\in\text{End}(\R^d)$ satisfies
$u=\lim_k\rho^{-n_k}\gamma_k$ with $\gamma_k\in\Gamma,$
$\rho^{-1}\leq\|u\|<1,$ then we have $|\det
u|=\lim_k\rho^{-n_kd}=0$ and $u(L_c)\subset
L_{ce^{i\alpha}}\cup\{0\}$ with $\alpha\in\R.$ From condition
$(H_2)$ we can choose $\gamma_k=\gamma^k\in\Gamma$ with
$W=\text{Ker}\,u\not=0$ having codimension 1. Since
$\text{Im}\,u=\R v_\gamma=\R a,$ $a\not= 0,$ we get that
$u(L_c)\subset H_\rho a\cup\{0\}.$ Since $u^{-1}(a)=b+W$ with
$b\in\R^d\setminus\{0\}$ we get that $L_c\subset W\cup
H_\rho(b+W).$ It follows that, in the quotient space $\R^d\slash
W,$ $L_c$ is projected on a set which is countable and
$H_\rho$-invariant. If $W_i,$ $1\leq i\leq r$ is a family of such
subspaces then the subspace $\bigcap_{i=1}^rW_i$ has the same
property. In fact $V\slash\bigcap_{i=1}^rW_i$ can be identified
with the diagonal subspace of $V\slash W_1\times\ldots\times
V\slash W_r,$ and so the projection of $L_c$ into
$V\slash\bigcap_{i=1}^rW_i$ can be identified with a subset of the
product of the projections of $L_c.$ Hence, such a projection is
countable and $H_\rho$-invariant. Since the intersection of any
family of subspaces with the above properties is a finite
intersection, there exists a minimum subspace $W_0$, which has the
considered properties. This subspace is unchanged when $c$ is
replaced by $ce^{i\alpha}$. As a consequence the condition
$\gamma(L_c)\subset L_{ce^{i\theta}}$ for $\gamma\in\Gamma$
implies that $\gamma(W_0)=W_0.$ Since $\pi(L_c)=L_\Gamma$ and
$L_\Gamma$ is uncountable (see Lemma \ref{4.3}) $W_0$ is proper.
This contradicts the irreducibility of $\Gamma,$ and so $t=0,$
$e^{-i\theta}\varphi(\gamma.x)=\varphi(x)$ for every $x\in
L_\Gamma$ and for every $\gamma\in S.$ Therefore, for every $n,$
and for all $\gamma_i\in S,$ and for every $x\in L_\Gamma$ we have
$e^{-in\theta}\varphi(\gamma_1\ldots\gamma_n.x)=\varphi(x).$ Since
$\Gamma=\bigcup_0^\infty S^n$ satisfies $(H_2),$ we can deduce
that $\varphi\equiv\text{const}$ on $L_\Gamma,$ and
$e^{i\theta}=1.$ In fact, suppose that there are two different
points $x$ and $y$ in $L_\Gamma$ such that
$\varphi(x)\not=\varphi(y).$ Then
\begin{equation*}
\varphi(x)=e^{-in\theta}\varphi(\gamma_1\ldots\gamma_n.x)\not=
e^{-in\theta}\varphi(\gamma_1\ldots\gamma_ny)=\varphi(y)
\end{equation*}
and so,
\begin{equation}\label{xy}
0<|\varphi(x)-\varphi(y)|=|\varphi(\gamma_1\ldots\gamma_n.x)-\varphi(\gamma_1\ldots\gamma_ny)|
\end{equation}
for every $n\in\N$ and all $\gamma_i\in S.$ By proximality of
$\Gamma$ on $\pr(V)$ (see Proposition \ref{prox onpr}) we have
that there exists a sequence $\{\gamma_i\}_1^\infty$ such that
\begin{equation*}
\lim_n\delta(\gamma_1\ldots\gamma_n.x,\gamma_1\ldots\gamma_n.y)=0.
\end{equation*}
By continuity of $\varphi$ we a get a contradiction to \eqref{xy}.
\end{proof}

The following corollary which is a complement to Proposition 3
p.~45 in \cite{GR} clarifies the role of the aperiodicity
hypotheses of $\Gamma$ considered by Kesten in \cite{K}, Guivarc'h
and Raugi in \cite{GR} (Proposition 3 and Lemma p. 45), Lalley in
\cite{L} (Corollaries 11.3, 11.4), Eberlein in \cite{Eb} and
Dal'bo in \cite{Dalbo} in the context of lengths of closed
geodesics in negative curvature. For an extension of these results
and their use in the more general setting of semisimple groups see
\cite{Benoist} and \cite{GS}. It explains why aperiodicity
conditions are not explicitly stated in Theorem \ref{main}, as in
\cite{B} and \cite{03b}.

\begin{corollary}\label{compcoll}
Suppose $\Gamma\subset GL(d,\R)$ is a sub-semigroup which
satisfies ($H_1$), ($H_2$) and denote
$S_\Gamma=\{\log|\lambda_g|:g\in\Delta_\Gamma\}.$ Then $S_\Gamma$
generates a dense subgroup of $\R.$
\end{corollary}
For the proof, we need three lemmas, the first of them being well
known (see \cite{Bow}, pp.~90--94).
\begin{lemma}\label{b}
Let $A$ be a finite set, $\Omega$ $\text{the compact metric space
$A^{\N}$}$ and $\theta$ the shift transformation on $\Omega$ given
by $(\theta\omega)_k=\omega_{k+1},$ $k\in\N.$ For a function
$\varphi$ on $\Omega$ we denote:
\begin{equation*}
S_n\varphi(\omega)=\sum_0^{n-1}\varphi\circ\theta^k(\omega).
\end{equation*}
Suppose $\varphi$ is H\"{o}lder continuous, and for any periodic
point $\omega$ of period $p,$ the sum $S_p\varphi(\omega)$ belongs
to $\Z.$ Then there exists a H\"{o}lder $\Z$-valued function
$\varphi^\prime$ on $\Omega$ and a H\"{o}lder function $\psi$ such
that $\varphi=\varphi^\prime+\psi-\psi\circ\theta.$
\end{lemma}
\begin{lemma}\label{c}
Suppose $g,h\in GL(d,\R)$ are such that $h$ is proximal and
$g.h^+\not\in V_h^<.$ Then, for $n=2p$ large, $gh^n$ is proximal
and
\begin{align*}
\lim_n(gh^n)^+&=g.h^+,\\
\lim_n V_{gh^n}^<&=V_h^<.
\end{align*}
\end{lemma}
\begin{proof}
We consider the sequence of linear maps $u_n=\frac{h^n}{\|h^n\|}$
and observe that $u_n$ converges towards a map $\pi_h$ which is
proportional to the projection on $\R v_h$ along the subspace
$V_h^<.$ Hence $gu_n$ converges towards $g\pi_h.$ We have
\begin{equation*}
\text{Im($g\pi_h$)}=\R g.h^+,\,\,\text{Ker($g\pi_h$)}=V_h^<.
\end{equation*}
Hence, if $g.h^+\not\in V_h^<,$ then $g\pi_h$ is collinear to a
projection onto a one-di\-men\-sio\-nal subspace. Since $g\pi_h$
has a simple dominant eigenvalue, the same is true for $gu_n$ for
$n$ large. Therefore, for $n$ large, $gu_n$ is proximal and
moreover, we have the required convergence.
\end{proof}
\begin{lemma}\label{d}
Suppose $\Gamma\subset GL(V)$ is a sub-semigroup and satisfies
conditions $(H_1),$ $(H_2).$ Then, there exists
$a,b\in\Delta_\Gamma$ such that $a^+\not=b^+,$ $V_a^<\not=V_b^<$
and $a^+\not\in V_b^<,$ $b^+\not\in V_a^<.$
\end{lemma}
\begin{proof}
We consider the transposed semigroup $\Gamma^t$ acting on the dual
space $V^*$ and the projective space $\pr(V^*).$ From Remark
\ref{requiv} iii) conditions $(H_1)$ and $(H_2)$ remain valid for
$\Gamma^t$ and we can consider the corresponding limit set
$L_{\Gamma^t}=L_\Gamma^*.$ We fix $a\in\Delta_\Gamma$ and we
observe that we can find $b\in\Delta_\Gamma$ such that
$V_b^<\not=V_a^<,$ $a^+\not\in V_b^<.$ Otherwise there will be a
dense subset of $L_\Gamma^*$ contained in a union of the two
projective subspaces defined by $V_a^<$ and $a^+$ in $\pr(V^*).$
Hence $L_\Gamma^*$ itself will be contained in such an union. But,
from Lemma \ref{4.3} below this is impossible. If $b^+\not\in
a^+\cup V_a^<,$ we have found the required pair $(a,b).$ If not,
we consider $g\in\Gamma$ and the sequence $gb^n,$ $n\in 2\N.$ In
view of Lemma \ref{l2.8} we can choose $g\in\Gamma$ such that:
\begin{equation*}
g.b^+\not\in V_b^<\cup V_a^<\cup a^+.
\end{equation*}
Then we can apply Lemma \ref{c} and replace $b$ by $gb^n=b^\prime$
for $n$ large. Under this condition, $V_{b^\prime}^<$ is close to
$V_b^<$ and the relations are still satisfied. Since
$(b^\prime)^+$ is close to $g.b^+$ and $g.b^+\not\in a^+\cup
V_a^<,$ the condition $(b^\prime)^+\not\in a^+\cup V_a^<$ is also
satisfied. Hence, we can take $(a,b^\prime)$ as the required pair.
\end{proof}
\begin{proof}[Proof of Corollary \ref{compcoll}.]
Since $\Gamma$ satisfies $(H_1)$ and $(H_2)$ we can choose
$a_1,a_2$ in $\Gamma$ according to the Lemma \ref{d}. Let
$C_1,C_2$ be two closed and disjoint neighborhoods of
$a_1^+,a_2^+$ in $\pr(V)$ such that $(C_1\cup
C_2)\cap(V_{a_1}^<\cup V_{a_2}^<)=\emptyset$ and let $o$ be a
point outside $V_{a_1}\cup V_{a_2}^<.$ Then, for $i=1,2$:
\begin{align*}
\lim_na_i^n.(C_1\cup C_2)=&a_i^+\\
\lim_na_i^n.o=&a_i^+.
\end{align*}
If we take $n$ large and set $a=a_1^n,$ $b=a_2^n,$ we have
\begin{equation}\label{3r}
a.o\in C_1,\,b.o\in C_2,\,a.(C_1\cup C_2)\subset \text{Int}\,
C_1,\,b.(C_1\cup C_2)\subset \text{Int}\,C_2.
\end{equation}
It follows from \eqref{3r} that the semigroup $\Gamma(a,b)$
generated by $a,b$ is free. In order to prove Corollary
\ref{compcoll} we can suppose $\Gamma=\Gamma(a,b).$ We consider
the trivial metric $\delta$ on $\{a,b\}$ and endow
$\Omega=\{a,b\}^\N$ with the metric
$\delta(\omega,\omega^\prime)=\sum_1^\infty
2^{-k}\delta(\omega_k,\omega_k^\prime).$ We define a homeomorphism
$\mathfrak{z}$ between $\Omega$ and $L_\Gamma$ as follows. We
observe that if $\omega=(a_k)_{k\in\N}$ and $a_k\in\{a,b\},$ then
it follows from \eqref{3r} that the sequence $a_1\ldots a_n.o$
converges to $\mathfrak{z}(\omega)\in C_1\cup C_2.$ It is easy to
verify that $\mathfrak z$ is a bi-H\"{o}lder homeomorphism, hence
we can transfer properties of $(\Omega,\theta)$ to the action of
$\Gamma$ on $L_\Gamma.$ We consider $\mathfrak z(\omega)$ as a
unit vector in $V$ and we observe that, by definition of
$\mathfrak z$:
\begin{equation*}
a_1(\omega)\mathfrak z(\theta\omega)=\|a_1(\omega)\mathfrak
z(\theta\omega)\|\mathfrak z(\omega).
\end{equation*}
It follows that, if we denote
\begin{align*}
\varphi(\omega)=&\log\|a_1(\omega)\mathfrak z(\theta\omega)\|,\\
S_n\varphi(\omega)=&\sum_0^{n-1}\varphi(\theta^k\omega)
\end{align*}
we have, with $\gamma=a_1\ldots a_{n-1}\in\Gamma$ and $x=\mathfrak
z(\theta^n\omega)\in L_\Gamma:$
\begin{equation*}
S_n\varphi(\omega)=\log\|\gamma x\|.
\end{equation*}
Given a H\"{o}lder function $\psi$ on $\Omega$ we define a
H\"{o}lder function $\bar\psi$ on $L_\Gamma$ by
$\bar{\psi}[\mathfrak z(\omega)]=\psi(\omega)$ and we have also
$\psi(\omega)-\psi(\theta^n\omega)=\bar{\psi}(\gamma.x)-\bar{\psi}(x).$
In particular, if $\omega\in\Omega$ is periodic with period $p,$
($\theta^p\omega=\omega$), then $\mathfrak z(\omega)$ is a
dominant eigenvector of $\gamma=a_1\ldots a_{p-1}$ and the
corresponding eigenvalue $\lambda_\gamma$ satisfies:
\begin{equation*}
\log|\lambda_\gamma|=S_p\varphi(\omega).
\end{equation*}
If $S_\Gamma$ do not generate a dense subgroup of $\R,$ then, for
some positive $c$ we have $S_\Gamma\subset c\Z,$ hence
$S_p\varphi(\omega)\in c\Z$ for any periodic point $\omega$ and we
can apply Lemma \ref{b} to the function $\frac{1}{c}\varphi.$ In
particular, the function $e^{2i\pi\varphi\slash c}$ can be written
in the form $e^{2i\pi(\psi-\psi\circ\theta)},$ where $\psi$ is a
H\"{o}lder function on $\Omega.$ We can define $\bar{\psi}$ on
$L_\Gamma$ as above and write $u(x)=e^{2i\pi\bar{\psi}}.$ Then $u$
is continuous and we obtain with $\gamma=a_1\ldots a_n,$
$x=\mathfrak z(\theta^n\omega)$:
\begin{equation*}
\|\gamma x\|^{2i\pi\slash c}=\frac{u(\gamma.x)}{u(x)}.
\end{equation*}
We extend  $u$ to $\pr(V)$ as a continuous function, again denoted
by $u.$ Then we have
\begin{equation*}
\forall\,x\in L_\Gamma,\,\forall\,\gamma\in\Gamma,\,\,\|\gamma
x\|^{2i\pi\slash c}=\frac{u(\gamma.x)}{u(x)}.
\end{equation*}
In view of the Proposition \ref{2.2}, this implies
$\frac{2i\pi}{c}=0,$ $u=1$ and this is impossible.
\end{proof}

\section{Random walks on a vector space and its factor spaces}\label{orbits}
In this section, relying strongly on \cite{F1}, \cite{GR} (see
also \cite{GR1}), we develop the random walk approach to the study
of $\Gamma$-orbits on $V\setminus\{0\}$ and other related
$\Gamma$-spaces. The main new results are Theorems \ref{5.1} and
\ref{tm} and their corollaries. They give Theorem
{\ref{essentialtool} and Theorem \ref{maintool} of the
Introduction. In particular, Corollary \ref{tm1} is one of the
main tools for the study of $\Gamma$-orbits on the torus $\T^d$
with $\Gamma\subset \m,$ i.e., for Theorem \ref{main}.

Let $\mu$ be a probability measure on $G=GL(V),$ $\Gamma_\mu$ the
closed sub-semigroup generated by the support $S_\mu$ of $\mu.$ We
denote by $M^1(X)$ the set of probability measures on a given
Polish space $X.$ We denote $\Omega=S_\mu^\N$ and we consider the
probability measure $\p_\mu=\mu^{\otimes\N}$ on $\Omega;$ the
shift $\theta$ on $\Omega$ given by
$(\theta\omega)_k=\omega_{k+1},$ ($k\in\N$) preserves $\p_\mu$ and
the components $\omega_k=g_k(\omega)$ of $\omega$ are i.i.d.
$G$-valued random variables of law $\mu.$ From Markov-Kakutani
theorem, there exists a probability measure $\nu$ on $\pr(V)$
which is $\mu$-stationary, i.e.,
\begin{equation*}
\mu*\nu=\int g.\nu d\mu(g)=\nu.
\end{equation*}
We are going to establish that $[\pr(V),\nu]$ is a $\mu$-boundary
(see \cite{F1}), i.e.,
\begin{equation*}
\lim_ng_1g_2\ldots g_n.\nu=\delta_{z_\omega},
\end{equation*}
where $z_\omega\in\pr(V).$ This will allow us to derive some
properties of the typical sequences:
\begin{equation*}
S_n=g_ng_{n-1}\ldots g_1\text{ and }X_n=g_1g_2\ldots g_n,
\end{equation*}
and transposed maps $S_n^t$ and $X_n^t.$
\begin{lemma}\label{4.3}
Assume that $\Gamma=\Gamma_\mu$ satisfies condition $(H_1),$ and
let $\nu$ be a $\mu$-stationary measure. Then $\nu$ gives zero
mass to every projective subspace. Furthermore, if $\Gamma$
satisfies also ($H_2$), then $L_\Gamma$ is not contained in a
countable union of subspaces.
\end{lemma}
\begin{proof}
Let $W$ be a projective subspace of minimal dimension such that
$\nu(W)>0.$

Define,
\begin{equation}\label{eta}
\eta=\eta_\mu=\sum_{k\geq 1}(1\slash 2^k)\mu^{*k}
\end{equation}
and consider the function $f(g)=g.\nu(W)=\nu(g^{-1}.W).$ This
function is $\mu$\textit{-harmonic}, i.e. satisfies
\begin{equation*}
\int f(gh)d\mu(h)=\int f(gh)d\eta(h)=f(g)
\end{equation*}
and reaches its maximum. In fact, the hypothesis on $W$ gives
$\nu(g.W\cap g^\prime.W)=0$ if $g.W\not=g^\prime.W,$ so the set of
$g.W$ such that $\nu(g.W)>\delta$ is finite for every $\delta.$
Then if $f(g_0)=\sup_{g\in G}f(g),$ the equation $f(g_0)=\int
f(g_0h)d\eta(h)$ gives $f(g_0h)=f(g_0),$ $\eta$-a.e. Consequently,
there exists a finite set $E$ of projective subspaces such that
$h^{-1}g_0^{-1}.W$ belongs to $E$ for a.e. $h,$ and therefore for
every $h\in\Gamma^{-1}.$ Then $\Gamma$ permutes a finite set of
projective subspaces and the strong irreducibility of $\Gamma$
gives the contradiction.

If ($H_2$) is also satisfied by $\Gamma,$ then $L_\Gamma$ is well
defined. From Markov-Kakutani theorem, we know that there exists a
$\mu$-stationary measure $\lambda$ such that
$\lambda(L_\Gamma)=1,$ hence $S_\lambda\subset L_\Gamma.$ Since
$\lambda$ gives zero measure to every subspace, the same is true
for a countable union of subspaces, hence $L_\Gamma$ cannot be
contained in such an union.
\end{proof}
\begin{proposition}\label{p}
Let $\nu$ be a $\mu$-stationary measure on $\pr(V)$ and $\eta$ be
as in \eqref{eta}. Then the sequence $g_1\ldots g_n.\nu$ converges
$\p_\mu$-a.e. and for $\p_\mu\otimes\eta$-a.e.
$(\omega,g)\in\Omega\times G$ we have
\begin{equation*}
\lim_ng_1(\omega)\ldots g_n(\omega).\nu=\lim_ng_1(\omega)\ldots
g_n(\omega)g.\nu.
\end{equation*}
\end{proposition}
\begin{proof}
For a continuous function $\varphi$, we set
$F_\varphi(g)=g.\nu(\phi)$ and we observe that the relation $\int
g.\nu d\mu(g)=\nu$ gives $\int F_\varphi(gh) d\mu(h)=F_\varphi(g)$
and consequently that $F_\varphi(X_n)=g_1\ldots g_n.\nu(\varphi)$
is a bounded martingale. This martingale converges and, letting
$\varphi$ vary in a dense countable part of $\pr(V),$ we obtain
the convergence of $g_1\ldots g_n.\nu.$ In order to obtain the
second claim, it suffices to show that
$F_\varphi(X_ng)-F_\varphi(X_n)$ converges to zero,
$\p_\mu\otimes\mu^{*r}$-a.e. for every $r\geq1.$ But:
\begin{equation*}
\int
|F_\varphi(X_ng)-F_\varphi(X_n)|^2d\mu^{*r}(g)d\p(\omega)=\mu^{*n+r}(F_\varphi^2)-\mu^{*n}(F_\varphi^2)
\end{equation*}
because of the relation $\int
F_\varphi(X_ng)d\mu^{*r}(g)=F_\varphi(X_n).$ One deduces that:
\begin{equation*}
\int\int\sum_{n=0}^p|F_\varphi(X_ng)-F_\varphi(X_n)|^2d\p(\omega)d\mu^{*r}(g)\leq
2r\|\varphi\|_\infty
\end{equation*}
for every $r\geq 1;$ this proves the convergence
$\p\otimes\mu^{*r}$-a.e. of the series
$\sum|F_\varphi(X_ng)-F_\varphi(X_n)|^2$ and consequently the
convergence of $F_\varphi(X_ng)-F_\varphi(X_n)$ to zero
\end{proof}
In the sequel we are going to use concepts introduced in
\cite{F1}. Therefore, we recall them briefly.

To every linear transformation of $\R^d$ is associated a
\textit{quasi-projective} transformation acting on the lines of
$\R^d$ not contained in the kernel of the transformation. So we
have maps of $\pr^{d-1}$ defined outside a projective subspace:
these maps are continuous outside the exceptional subspace and are
limits, outside this subspace, of a sequence of projective
transformations. Furthermore, from every sequence of projective
transformations, we can extract a subsequence converging to a
quasi-projective one, outside a projective subspace.
\begin{theorem}\label{4.2}
Let $\nu$ be a $\mu$-stationary measure on $\pr(V).$ Assume that
$\Gamma=\Gamma_\mu$ satisfies conditions $(H_1)$ and $(H_2).$ Then
we have $\p_\mu$-a.e.
\begin{equation*}
\lim_n g_1g_2\ldots g_n.\nu=\delta_{z_\omega}.
\end{equation*}
In particular $\nu$ is unique and its support is $L_\Gamma.$
\end{theorem}
\begin{proof}
The proof goes like in \cite{GR}. For a fixed $\omega,$ we
consider the relation given by Proposition \ref{p},
\begin{equation*}
\theta(\omega)=\lim_ng_1(\omega)\ldots g_n(\omega).\nu=\lim_n
g_1(\omega)\ldots g_n(\omega)g.\nu,
\end{equation*}
which is true for $\p_\mu\otimes\eta$-a.e. $(\omega,g).$ One can
extract from $ g_1(\omega)\ldots g_n(\omega)$ a subsequence
converging outside a projective subspaces to a quasi-projective
map $\tau(\omega).$ As $\nu$ gives zero-measure to any projective
subspace (Lemma \ref{4.3}), one has from Proposition \ref{p}
above: $\tau(\omega).\nu=\tau(\omega)g.\nu=\theta(\omega),$ for
$\eta$-a.e. $g$, and therefore for all $g\in\Gamma.$ As $\Gamma$
satisfies $(H_1)$ and $(H_2),$ one can find a sequence
$g_n\in\Gamma$ such that $g_n.\nu$ converges to a Dirac measure
$\delta_z$ with $z$ belonging to the open set of continuity of
$\tau(\omega).$ Then, in the limit
$\theta(\omega)=\tau(\omega).\delta_{z}.$ This proves that
$\theta(\omega)$ is a Dirac measure $\delta_{z_\omega}.$ The law
of random variable $z$ is necessarily $\nu$ by the martingale
convergence theorem. Since $z$ is independent of the choice of the
$\mu$-stationary measure $\nu$ we get the uniqueness of $\nu.$

Clearly, $S_\nu$ is closed and $\Gamma$-invariant. Hence,
Proposition \ref{4.1} gives $S_\nu\supset L_\Gamma.$ The
Markov-Kakutani theorem and uniqueness of $\nu$ give, as in the
proof of Lemma \ref{4.3}, $\nu(L_\Gamma)=1,$ hence
$S_\nu=L_\Gamma.$
\end{proof}
\begin{corollary}\label{c4}
Let $\rho$ ($\rho^*$ resp.) be a probability measure on $\pr(V)$
($\pr(V^*)$ resp.) which gives zero mass to every projective
subspace. Then we have $\p_{\mu}$-a.e.
\begin{equation*}
\lim_n g_1\ldots g_n.\rho=\delta_{z_\omega}\text{
($\lim_ng_1^t\ldots g_n^t.\rho^*=\delta_{z_\omega^*}$ resp.)}
\end{equation*}
In particular:
\begin{equation*}
\lim_n g_1\ldots g_n.m =\delta_{z_\omega}\text{
($\lim_ng_1^t\ldots g_n^t.m^*=\delta_{z_\omega^*}$ resp.),}
\end{equation*}
where $m$ ($m^*$ resp.) is the $K$-invariant probability measure
on $\pr(V)$ ($\pr(V^*)$ resp.)
\end{corollary}
\begin{proof}
We observe that $\frac{g_1\ldots g_n}{\|g_1\ldots
g_n\|}\in\text{End}(V)$ has norm one. We consider an arbitrary
convergent subsequence:
\begin{equation*}
u=\lim_k\frac{g_1\ldots g_{n_k}}{\|g_1\ldots g_{n_k}\|}.
\end{equation*}
Clearly $u\not=0,$ since $\|u\|=1.$ We note that $u$ defines a
continuous map from $\pr(V)\setminus\text{Ker}\, u$ onto $\pr(V).$
We will denote it again by $u,$ and observe that, since
$\rho(\text{Ker}\, u)=0,$ $u.\rho$ is well defined, and from
dominated convergence:
\begin{equation*}
u.\rho=\lim_k g_1\ldots g_{n_k}.\rho.
\end{equation*}
In particular from Theorem \ref{4.2} and Lemma \ref{4.3}:
$u.\nu=\delta_{z_\omega}.$ This means that the linear map $u$ has
rank one and satisfies $u(\pr(V)\setminus\text{Ker}\,
u)=\delta_{z_\omega}.$ Hence $u.\rho=\delta_{z_\omega}.$ The
convergent subsequence chosen above was arbitrary, hence:
\begin{equation*}
\lim_ng_1\ldots g_n.\rho=\delta_{z_\omega}.
\end{equation*}
In particular, we have the above convergence for $\rho=m.$

The results for $\pr(V^*),$ $\rho^*,$ $m^*,$ $z_\omega^*$ follow
from $\Gamma_{\mu^*}=(\Gamma_\mu)^t$ and Remark \ref{requiv} iii).
\end{proof}
Recall $m$ is the $K$-invariant probability measure on $\pr(V),$
where $K=SO(d,\R).$ One says that a sequence $f_n\in GL(d,\R)$ has
the \textit{contraction property on $\pr(V)$ towards $z$} if the
sequence of measures $f_n.m$ on $\pr(V)$ converges weakly towards
$\delta_z.$ A point $z\in\pr(V)$ will be identified with a vector
of norm one, defined up to a sign.

We will use, as in \cite{GR} and \cite{GR1}, the
$K\overline{A^+}K$ decomposition of $g\in GL(d,\R),$
$g=kak^\prime,$ $k$ and $k^\prime$ are orthogonal matrices and
$\overline{A^+}\ni a=\text{diag}(a^1,\ldots,a^d)$ with $a^1\geq
a^2\geq\ldots\geq a^d$ and $(e_1,\ldots,e_d)$ denotes the
canonical basis of $\R^d.$ In particular, if $g\in SL(d,\R)$ then
$k,k^\prime\in K=SO(d,\R)$ and $a\in
\overline{A^+}=\{\text{diag}(a^1,\ldots,a^d):a^1\geq
a^2\geq\ldots\geq a^d>0\text{ and }\prod_1^da^i=1\}.$

If one writes the polar decomposition of $f_n$ as
$f_n=k_na_nk_n^\prime,$ where $k_n,k_n^\prime\in K,$
$a\in\overline{A^+},$ one sees that the contraction property is
equivalent to $a_n^{i}=o(a_n^{1}),$  $1<i\leq d,$ $\lim_n k_n.
e_1=z.$

In proposition below and its corollary, the point $z\in\pr(V)$ is
considered as a unit vector, hence $|z(x)|$ is well defined for
$x\in V.$
\begin{proposition}\label{4.5}
Assume that $f_n\in GL(V)$ is a sequence such that $f_n^{t}$ has
the contraction property on $\pr(V^*)$ towards $z\in\pr(V^*).$
Then for any $x,y\in\pr (V),$ one has the following convergence
\begin{equation*}
\lim_n\frac{\|f_n(x)\|}{\|f_n\|}=|z(x)|,
\end{equation*}
\begin{equation*}
\lim_nz(x)z(y)\frac{\delta(f_n.x,f_n.y)}{\delta(x,y)}=0.
\end{equation*}
The second convergence is uniform when $x,y$ belong to a compact
subset of $\pr(V)\setminus\text{Ker $z.$}$

If $f_n\in SL(V),$ then one has for every $x\not\in\text{Ker
$z,$}$ $\lim_n\|f_n(x)\|=+\infty.$
\end{proposition}
\begin{proof}
Recall that the distance between elements $\bar{u}=\pi(u)$ and
$\bar{v}=\pi(v)$ in $\pr(V)$ is equal to
$\delta(\bar{u},\bar{v})=\|u\wedge v\|\slash \|u\|\|v\|.$

One writes $f_n=k_na_nk_n^\prime$ as above, with
$k_n,k_n^\prime\in K$ and $a_n\in\overline{A^+}.$ From the
hypotheses we get:
\begin{equation*}
\lim_n k_n^{\prime -1}.\bar{e_1}=z,\, a_n^{i}=o(a_n^1),\, (i>1).
\end{equation*}
Writing $x=\sum_1^d x_ie_i,$ we get:
\begin{equation*}
\|f_nx\|^2=\|a_nk_n^\prime x\|^2=\sum_1^d(a_n^i)^2|\langle
k_n^\prime x,e_i\rangle |^2\geq (a_n^1)^2|\langle k_n^{\prime}
x,e_1\rangle|^2.
\end{equation*}
Since the norm of $f_n$ is $a_n^1,$ we get
\begin{equation*}
\begin{split}
\lim_n\frac{\|f_nx\|^2}{\|f_n\|^2}&=\lim_n|\langle x,k_n^{\prime
-1}e_1\rangle|^2+\lim_n\sum_{i>1}\left(\frac{a_n^i}{a_n^1}\right)^2|\langle
k_n^\prime x, e_1\rangle|^2\\&=\lim_n|\langle x,k_n^{\prime
-1}e_1\rangle|^2=|z(x)|^2.
\end{split}
\end{equation*}
Also
\begin{equation*}
\begin{split}
\|f_n(x)\wedge f_n(y)\|^2=&\sum_{i\not=j}(a_m^ia_n^j)^2|\langle
k_n^\prime(x\wedge y), e_i\wedge e_j\rangle|^2\\
\|f_n(x)\wedge f_n(y)\|\leq & d^2a_n^1a_n^2\|x\wedge y\|.
\end{split}
\end{equation*}
Therefore,
\begin{equation*}
\frac{\delta(f_n.x,f_n.y)}{\delta(x,y)} =\frac{\|f_nx\wedge
f_ny\|}{\|f_nx\|\|f_ny\|\|x\wedge y\|} \leq
d^2\frac{a_n^2}{a_n^1}\frac{1}{|\langle k_n^{\prime }x,
e_1\rangle\langle k_n^{\prime }y, e_1\rangle|}
\end{equation*}
and
\begin{equation*}
\lim_n|z(x)||z(y)|\frac{\delta(f_n.x,f_n.y)}{\delta(x,y)}=0.
\end{equation*}
The uniformity of the required convergence is clear from the
previous formula, if $z(x)z(y)\not\equiv 0.$

In order to obtain the last assertion, it suffices to show, in
view of the first statement, that $\|f_n\|$ converges to
$+\infty.$ The relation $a_n^2=o(a_n^1)$ implies $\det
f_n=\prod_{i=1}^da_n^i=o(a_n^1).$ Since $\det f_n=1,$ we get
$\lim_n\|f_n\|=\lim_n a_n^1=+\infty.$
\end{proof}
\begin{corollary}[see \cite{GR,GR1}]\label{4.6}
If $\mu,$ $z_\omega^*$ are as in Theorem \ref{4.2} and Corollary
\ref{c4}, then, as $n$ tends to infinity, we have uniformly in
$x,y\in\pr(V):$
\begin{equation*}
\lim_n\e_\mu\delta(S_n.x,S_n.y)=0
\end{equation*}
and $\p_\mu$-a.e.
\begin{equation*}
\lim_n\frac{\|S_nx\|}{\|S_n\|}=|z_\omega^*(x)|.
\end{equation*}
If $\mu\in M^1(SL(V)),$ one has, for every $x\in V$ and
$\p_\mu$-a.e. $\lim_n\|S_nx\|=+\infty.$
\end{corollary}
\begin{proof}
Note that Corollary \ref{c4} implies that $S_n^t(\omega)$ has the
contraction property towards $z_\omega^*.$ Hence Proposition
\ref{4.5} implies:
\begin{equation*}
\lim_n\frac{\|S_nx\|}{\|S_n\|}=|z^*_\omega(x)|.
\end{equation*}

For the first convergence it suffices to show that for any
sequence $x_n,y_n\in\pr(V)$ we have $\p_\mu$-a.e.
\begin{equation}\label{number}
\lim_n\delta(S_n.x_n,S_n.y_n)=0.
\end{equation}
One can suppose that $\lim_nx_n=x$ and $\lim_ny_n=y.$ From
Corollary \ref{c4} and Lemma \ref{4.3}, one knows that the law of
$z_\omega^*=\lim_nS_n^t.m^*$ gives zero measure to every subspace.
Hence, for almost every $\omega\in\Omega,$ $x$ and $y$ are not in
$\text{Ker}\,z_\omega^*,$ and the same is true for $x_n,y_n$ for
large $n.$ Then \eqref{number} follows from dominated convergence.

The last assertion is proved as follows. If $x$ is fixed, then
$\p_\mu$-a.e. as above $|z_\omega^*(x)|\not=0,$ hence from
Proposition \ref{4.5}, $\lim_n\|S_nx\|=+\infty.$
\end{proof}

Now we are going to study stationary measures on factor spaces of
$V\setminus\{0\}.$

Let $c>1 $ be fixed and let us denote by $\pr_c(V)=\pr_c^{d-1}$
the factor space of $V\setminus\{0\}$ by the multiplicative
subgroup of $\R^*$:
\begin{equation*}
\pm c^\Z:=\{\pm c^n:n\in\Z\},
\end{equation*}
and by $\T_c$ the 1-torus $\T_c=\R^*\slash\pm c^\Z.$

We can consider the projection from $V\setminus\{0\}$ to $
\pr(V)\times\T_c$ given by
\begin{equation*}
v\mapsto (\bar{v},\|v\|^{i\alpha}),
\end{equation*}
where $\alpha=2\pi\slash\log c$ and we observe that $\pr_c(V)$ is
then naturally identified with $\pr(V)\times\T_c.$ Hence a point
of $\pr_c(V)$ will be written as $v=(\bar{v},z),$ where
$\bar{v}\in\pr(V)$ is the projection of $v$ and
$z=\|v\|^{i\alpha}.$ The action of $g\in G=GL(V)$ on $\pr_c(V)$
can then be written as
\begin{equation*}
g.v=g.(\bar{v},z)=(g.\bar{v},z\|gv\|^{i\alpha}).
\end{equation*}
$\R^*$ acts also on this space and the two actions commute. The
corresponding formula is
\begin{equation*}
t.(\bar{v},z)=(\bar{v},z|t|^{i\alpha}),\;t\in\R^*.
\end{equation*}

We denote by $\lambda_c=dz$ the normalized Lebesgue measure on
$\T_c$ and observe that every measure of the form
$\nu\otimes\lambda_c,$ where $\nu\in M^1(\pr(V))$ is invariant
under the action of $\R^*$ on $\pr_c(V).$ Furthermore, if $\mu\in
M^1(G)$ and $\nu\in M^1(\pr(V))$ is $\mu$-stationary, then
$\nu\otimes\lambda_c$ is also $\mu$-stationary. If
$L_\Gamma\subset\pr(V)$ is the limit set of $\Gamma$
($\Gamma=\Gamma_{\mu}$ satisfies $(H_1)$ and $(H_2)$ conditions),
then $L_\Gamma(c)=L_\Gamma\times\T_c$ is a closed and
$\Gamma$-invariant subset of $\pr(V)\times \T_c.$
\begin{theorem}\label{5.1}
Assume that $\mu\in M^1(G)$ is such that $\Gamma=\Gamma_{\mu}$
satisfies conditions $(H_1)$ and $(H_2)$. Then, with the above
notations, for every $\psi\in C[\pr_c(V)]$ the sequence
$\check{\mu}^{*n}*\psi$ converges uniformly to
$(\nu\otimes\lambda_c)(\psi),$ where $\nu$ is $\mu$-stationary
measure on $\pr(V).$ Furthermore, for any $v\in\pr_c(V)$ we have
the following a.e. convergence:
\begin{equation*}
\lim_n\delta^c(S_n(\omega).v,L_\Gamma(c))=0,
\end{equation*}
where $\delta^c$ is the distance on $\pr(V)\times\T_c$ given by
\begin{equation}\label{dispc}
\delta^c(v,v^\prime)=\delta(\bar{v},\bar{v^\prime})+|z-z^\prime|,
\end{equation}
and $v=(\bar{v},z),$ $v^\prime=(\bar{v^\prime},z^\prime).$
\end{theorem}
\begin{corollary}\label{c5.3}
Assume $\Gamma\subset G$ is a sub-semigroup of $G$ which satisfies
property $(H_1)$ and $(H_2)$ and $c>1$ is fixed. Then the closed
$\Gamma$-invariant subset $L_\Gamma(c)=L_\Gamma\times \T_c$ of
$\pr_c(V)$ is the unique minimal set. Furthermore, any $\mu\in
M^1(G)$ such that $\Gamma=\Gamma_{\mu}$ satisfies conditions
$(H_1)$ and $(H_2)$ has a unique stationary measure on $\pr_c(V).$
\end{corollary}
Clearly, Theorem \ref{5.1} and its corollary imply Theorem
\ref{essentialtool} of the Introduction.

For the proof of Theorem \ref{5.1} we need three lemmas.
\begin{lemma}\label{lemma5.4}
If $\mu$ is as in Theorem \ref{5.1} then for $x,y\in V,$
\begin{equation*}
\lim_{y\to
x}\limsup_n\e_\mu\left|\left(\frac{\|S_nx\|}{\|S_ny\|}\right)^{i\alpha}-1\right|=0.
\end{equation*}
\end{lemma}
\begin{proof}
>From Corollary \ref{4.6} we know that if $x_n\to x$ and $y_n\to
y,$ then:
\begin{equation*}
\lim_n\frac{\|S_nx_n\|}{\|S_ny_n\|}=\frac{|z_\omega(x)|}{|z_\omega(y)|}.
\end{equation*}
Hence, from dominated convergence:
\begin{equation*}
\limsup_n\e_\mu\left|\left(\frac{\|S_nx_n\|}{\|S_ny_n\|}\right)^{i\alpha}-1\right|=
\e_\mu\left|\left|\frac{z_\omega(x)}{z_\omega(y)}\right|^{i\alpha}-1\right|.
\end{equation*}
The formula in the lemma corresponds to the special case $x_n=x,$
$y=x.$
\end{proof}
\begin{lemma}\label{lemma5.5}
If $\mu$ is as in Theorem \ref{5.1}, for every $\psi\in
C(\pr_c(V))$ the sequence of functions $\check\mu^{*k}*\psi$ is
uniformly equicontinuous.
\end{lemma}
\begin{proof} One considers the distance $\delta^c$ on $\pr_c(V)$
given by \eqref{dispc}. Then, in view of the form of the action of
$G$ on $\pr_c(V)$:
\begin{equation*}
\delta^c(S_n.v,S_n.v^\prime)=\delta(S_n.\bar{v},S_n.\bar{v^\prime})
+\left|\left(\frac{\|S_n\bar{v}\|}{\|S_n\bar{v^\prime}\|}\right)^{i\alpha}-1\right|+|z-z^\prime|.
\end{equation*}
>From Lemma \ref{lemma5.4}, we get:
\begin{equation*}
\lim_n\delta^c(S_n.v,S_n.v^\prime)=\left|\left|\frac{\langle\bar{v},z_\omega^*\rangle}
{\langle\bar{v^\prime},z_\omega^*\rangle}\right|^{i\alpha}-1\right|+|z-z^\prime|.
\end{equation*}
Hence, using dominated convergence
\begin{equation}\label{q}
\limsup_n\e_\mu\delta^c(S_n.v,S_n.v^\prime)=|z-z^\prime|+\e_\mu\left(\left|\frac{\langle\bar{v},z_\omega^*\rangle}
{\langle\bar{v^\prime},z_\omega^*\rangle}\right|^{i\alpha}-1\right).
\end{equation}
The right hand side of this formula is uniformly small when
$\delta^c(v,v^\prime)$ is small.

Now, if $\psi\in C(\pr_c(V))$ is Lipschitz, with coefficient
$[\psi]:$
\begin{equation*}
|\check\mu^{*n}*\psi(v)-\check\mu^{*n}*\psi(v^\prime)|\leq\e_\mu[\delta^c[(S_n.v,S_n.v^\prime)][\psi].
\end{equation*}
Since Lipschitz functions are dense in $C(\pr_c(V))$ the above
inequality and \eqref{q} imply equicontinuity of the sequence
$\check\mu^{*n}*\psi(v)$ for $\psi\in C(\pr_c(V)).$
\end{proof}
\begin{lemma}\label{nowy lemat}
Suppose $\theta\in\R,$ $\eta\in C(\pr(V))$ and $\eta\not\equiv 0$
and  satisfies the equation
\begin{equation}\label{r1gt}
\int\eta(g.\bar{v})\|g\bar{v}\|^{i\alpha}d\mu(g)=e^{i\theta}\eta(\bar{v}).
\end{equation}
Then $\alpha=0,$ $\theta=0$ and $\eta=\text{const}$ on $\pr(V).$
\end{lemma}
\begin{proof}
Passing to absolute values in \eqref{r1gt} we get
\begin{equation}\label{r2gt}
|\eta(\bar{v})|\leq\int|\eta(g.\bar{v})|d\mu(g).
\end{equation}
Let $M=\{\bar{v}\in\pr(V):|\eta(\bar{v})|=\|\eta\|_\infty\}.$ Then
from \eqref{r2gt} the condition $\bar{v}\in M$ implies
$g.\bar{v}\in M$ $\mu$ a.e. Hence from continuity of $|\eta|,$ we
have $gM\subset M$ for every $g\in S_\mu$ and $\Gamma_\mu M\subset
M.$ Since $L_\Gamma$ is minimal we get $L_\Gamma\subset M.$ In
particular  for every $\bar{v}\in L_\Gamma,$
$|\eta(\bar{v})|=\|\eta\|_\infty.$ From strong convexity of the
unit disc in $\C$ we get from \eqref{r1gt} that
\begin{equation*}
\forall\,\bar{v}\in L_\Gamma,\,\forall g\in S_\mu,\,\,
\eta(g.\bar{v})\|gv\|^{i\alpha}=e^{i\theta}\eta(\bar{v}).
\end{equation*}
>From Proposition \eqref{2.2} it follows that $\alpha=0,$
$\theta=0$ and $\eta=\text{const}$ on $L_\Gamma.$

Now we have on $\pr(V):$
\begin{equation*}
\int\eta(g.\bar{v})d\mu(g)=\eta(\bar{v}).
\end{equation*}
We can suppose $\eta$ to be real and consider the set $M^\prime$
($M^{\prime\prime}$ resp.) of points where $\eta$ attains its
maximum (minimum resp.). As above we obtain $M^\prime\supset
L_\Gamma.$ Replacing $\eta$ by $-\eta$, we obtain also that
$M^{\prime\prime}\supset L_\Gamma,$ hence
$M^\prime=M^{\prime\prime}.$ We conclude that
\begin{equation*}
\forall\,\bar{v}\in\pr(V),\,\,\eta(\bar{v})=\text{const}.
\end{equation*}
\end{proof}
\begin{proof}[Proof of Theorem \ref{5.1}.]
We use the following result of \cite{Ro}. Let $P$ be a Markov
operator on the compact metric space $X,$ which preserves $C(X)$
and is equicontinuous, i.e., for any $\psi\in C(X),$ the sequence
$P^k\psi,$ $k\in\N$ is e\-qui\-con\-ti\-nuous. Then if $1$ is the
only eigenvalue of modulus one in $C(X),$ the sequence $P^k\psi$
converges uniformly. Here we have $P(x,\cdot)=\mu*\delta_x,$ and
$X=\pr_c(V).$ From Lemma \ref{lemma5.5} we know that $P$ is
equicontinuous. Suppose that $\eta\in C(X),$ $\eta\not\equiv 0$
satisfies $P\eta=e^{i\theta}\eta,$ i.e.,
\begin{equation*}
\int\eta(g.v)d\mu(g)=e^{i\theta}\eta(v)
\end{equation*}
for any $v$ in $\pr_c(V).$ Now we can consider the Fourier
coefficients ($k\in\Z$),
\begin{equation*}
\eta_k(\bar{v})=\int\eta(\bar{v},z)z^kd\lambda_c(z)
\end{equation*}
and we obtain
\begin{equation*}
\int\eta_k(g.\bar{v})\|g\bar{v}\|^{ik\alpha}d\mu(g)=e^{i\theta}\eta_k(\bar{v}).
\end{equation*}
>From Lemma \ref{nowy lemat} we get $e^{i\theta}=1,$
$\eta_k(\bar{v})=0$ for $k\not=0,$
$\eta_0(\bar{v})\equiv\text{const}.$ Hence $\eta=\text{const}$ on
$\pr_c(V).$ Now the result of \cite{Ro} give the uniform
convergence of the sequence $\psi_n=\check{\mu}^{*n}*\psi.$

Clearly, if $\lim_{n}\psi_{n}=\eta,$ one has
$P\eta=\check\mu*\eta=\eta$ and $\eta$ is continuous. From the
above result, we deduce $\eta\equiv\text{const}.$ Furthermore,
\begin{equation*}
\eta=(\nu\otimes\lambda_c)(\eta)=\lim_n(\nu\otimes\lambda_c)(\psi_{n})=(\nu\otimes\lambda_c)(\psi).
\end{equation*}
Hence the formula $\eta=(\nu\otimes\lambda_c)(\psi)$ and the
required convergence.

In order to prove the second statement of the theorem, notice that
since $L_\Gamma(c)$ is the inverse image of $L_\Gamma$ in
$\pr_c(V)$ we have:
\begin{equation*}
\delta^c(S_n(\omega).v,L_\Gamma(c))=\delta(S_n(\omega).\bar{v},L_\Gamma).
\end{equation*}
Proposition \ref{4.5} implies that, given $\bar{v}$ and $\bar{w}$
in $\pr(V),$ we have the a.e. convergence of the sequence
$\delta(S_n(\omega).\bar{v},S_n(\omega).\bar{w})$ to zero. If we
choose $\bar{w}$ in $L_\Gamma,$ then $S_n(\omega).\bar{w}\in
L_\Gamma,$ hence:
\begin{equation*}
\delta(S_n(\omega).\bar{v},L_\Gamma)\leq\delta(S_n(\omega).\bar{v},S_n(\omega).\bar{w}).
\end{equation*}
It follows:
$\lim_n\delta(S_n(\omega).\bar{v},L_\Gamma)=\delta(S_n(\omega).\bar{v},S_n(\omega).\bar{w})=0.$
\end{proof}
\begin{proof}[Proof of Corollary \ref{c5.3}.]
Suppose $\xi\in M^1(\pr_c(V))$ is another $\mu$-sta\-tio\-na\-ry
measure. Since $\psi_n=\check\mu^{*n}*\psi$ converges uniformly to
$(\nu\otimes\lambda_c)(\psi),$ we get
\begin{equation*}
\xi(\lim_n\psi_n)=\lim_n\xi(\check\mu^{*n}*\psi)=\xi(\psi).
\end{equation*}
Hence, $(\nu\otimes\lambda_c)(\psi)=\xi(\psi),$
$\nu\otimes\lambda_c=\xi$ and the uniqueness follows.

Suppose $\Delta$ is a closed $\Gamma_\mu$-invariant subset of
$\pr_c(V).$ Then from the Mar\-kov-Kakutani theorem, there is a
$\mu$-stationary measure carried by $\Delta$. From the uniqueness
of the stationary measure we get
\begin{equation*}
\Delta\supset\text{supp }\nu\otimes\lambda_c=L_\Gamma(c).
\end{equation*}
\end{proof}
\begin{theorem}\label{tm}
Suppose that $\Gamma$ is a sub-semigroup of $GL(d,\R)$ satisfying
conditions $(H_0),$ $(H_1)$ and $(H_2)$ and let $\Sigma$ be
$\Gamma$-invariant subset of $\tilde V\setminus\{0\}$ such that
$0$ is a limit point of $\Sigma.$ Then
\begin{equation}\label{v301}
\overline{\Sigma}\supset \tilde L_\Gamma\slash\{\text{Id
},\sigma\}.
\end{equation}
\end{theorem}
\begin{proof}
We denote by $\Sigma^\prime$ the inverse image of $\Sigma$ in
$V\setminus\{0\}.$ Let $u_0$ be a $\gamma$-dominant vector as in
Proposition \ref{rgs}, that is satisfying
\begin{equation}\label{v30}
\gamma^{\Z}u_0:=\{\gamma^ku_0:k\in\Z\}\subset\overline{\Sigma^\prime.}
\end{equation}
Applying Corollary \ref{c5.3} with $c=\lambda,$ where $\lambda$ is
the unique eigenvalue of $\gamma$ of maximum modulus, greater than
1 since $\gamma$ is expanding, we get that if $\bar{u_0}$ denotes
the projection of $u_0$ on $\pr_c(V)$ then
$\overline{\Gamma\bar{u_0}}\supset L_\Gamma(c).$ It follows that
if $\bar{y}\in L_\Gamma(c)$ is given, then there is a sequence
$\{\gamma_n\}\subset\Gamma,$ such that $\gamma_n.\bar{u_0}$
converges to $\bar{y}.$ This implies that there is a sequence of
integers $\{p_n\}$ such that $\lambda^{p_n}\gamma_n u_0\to y\in
V\setminus\{0\}$ but this implies
\begin{equation*}
\gamma_n\lambda^{p_n}u_0=\gamma_n\gamma^{p_n}u_0\to y.
\end{equation*}
But $\gamma^{p_n}u_0\in\overline{\Sigma^\prime}$ by \eqref{v30}.
Thus $y\in\overline{\Sigma^\prime}.$ Since $\bar{y}$ was an
arbitrary point from $L_\Gamma(c)$ we get that $\tilde
L_\Gamma\subset\overline{\Sigma^\prime}$ and \eqref{v301} is
proved.
\end{proof}
Clearly, Theorem \ref{tm} gives Theorem \ref{maintool} of the
Introduction.

Theorem \ref{tm} will be used below in the special case
$\Gamma\subset \m.$
\begin{corollary}\label{tm1}
Let $\Gamma$ be a sub-semigroup of $\m$ satisfying $(H_0),$
$(H_1)$ and $(H_2).$ Let $\Sigma$ be a $\Gamma$-invariant subset
of $\tilde V\setminus\{0\}$ such that $0$ is a limit point of
$\Sigma.$ Then $\overline{\Sigma}\supset \tilde
L_\Gamma\slash\{\text{Id, $\sigma$}\}.$
\end{corollary}

Theorem \ref{tm} does not give information on a general
$\Gamma$-orbit closure $\overline{\Gamma v},$ $v\in\tilde
V\setminus\{0\}$ if $0$ is not a limit point. On the other hand
Theorem \ref{5.1} and its corollary describe the behavior of a
general $\Gamma$-orbit in $\pr_c(V).$ Using more precise
informations on products of random matrices, i.e the renewal
theorem as in \cite{GS} (see also \cite{K}), one can go further
and describe the behavior at infinity of a general orbit $\Gamma
v\subset\tilde V\setminus\{0\}$ as follows. For any $c,d$ ($1\leq
c<d$) we denote by $\tilde V_{[c,d]}\subset\tilde V\setminus\{0\}$
the "$c$-shell" $\pr^{d-1}\times[c,d],$ by $\tilde
L_{\Gamma,c}\subset\tilde V_c:=\pr^{d-1}\times[1,c]$ the closed
subset $L_\Gamma\times[1,c].$ Then by the methods of \cite{CG01}
we can obtain the following
\begin{theorem}\label{ot}
Assume that the semigroup $\Gamma\subset GL(d,\R),$ $d>1,$
satisfies $(H_0),$ $(H_1)$ and $(H_2).$ Then, with the above
notations, for any $c>1,$ $v\in\tilde V\setminus\{0\}$ we have the
following convergence
\begin{equation*}
\lim_{t\to\infty}c^{-t}(\overline{\Gamma v}\cap\tilde
V_{[c^t,c^{t+1}]})=\tilde L_{\Gamma,c}.
\end{equation*}
\end{theorem}

This can be interpreted as "thickness" at infinity in the
direction of $\tilde L_\Gamma$ of the orbit closure
$\overline{\Gamma v}\subset\tilde V.$

Theorems \ref{tm} and Corollary \ref{tm1} can also be deduced from
Theorem \ref{ot}.
\begin{remark}
\emph{The conclusions in statements \ref{tm} to \ref{ot} are valid
also if $d=1,$ if one supposes the semigroup $\Gamma$ of $\R^*$ to
be non-lacunary. The corresponding aperiodicity condition in the
statements above is automatically satisfied if $d>1,$ because of
Corollary \ref{compcoll}.}
\end{remark}

\section{Proof of Theorem \ref{main}}\label{proof}
In order to prove the theorem, we use ideas of \cite{F} and
\cite{B}. The first step is to prove that if $\Sigma\subset\T^d$
is a closed $\Gamma$-invariant subset that contains $0\in\T^d$ as
a limit point, then $\Sigma=\T^d.$ Here we apply Corollary
\ref{tm1} to the inverse image $p^{-1}(\Sigma)$ of $\Sigma$ in
$\R^d.$

In the general case, we suppose $\Sigma$ to be infinite and we
construct other closely related closed $\Gamma$-invariant subsets
of $\T^d$ which contains $0.$ Then we use the special case to get
information on $\Sigma$ and we conclude that $\Sigma=\T^d.$
\subsection{The case when $0$ is a limit point of
$\Sigma$}\label{subs1}

The statement $\Sigma=\T^d$ will hold by Corollary \ref{tm1}
applied to $p^{-1}(\Sigma)$ if we are able to see that $\tilde
L_\Gamma$ contains at least one ray which is not contained in a
rational subspace. But the set of rational subspaces is countable
and, by Lemma \ref{4.3} $L_\Gamma$ is not contained in a countable
union of subspaces. The result follows.

We can observe that the set $\tilde L_\Gamma$ is very large, since
it was proved in \cite{CG00} that $L_\Gamma$ has strictly positive
Hausdorff dimension.

\subsection{The general case}\label{subs2}
In order to show that the above case is the only one, we make use
of previous ideas from \cite{F} and \cite{B}.

If $\gamma\in \m$ and $m\in\N$ is fixed we write
\begin{equation*}
\gamma\equiv\text{Id}\text{ (mod $m$) $\iff$
$\gamma-\text{Id}=mA,$}
\end{equation*}
with $A\in M(d,\Z):=\{d\times d\text{ matrices with integer
entries}\}.$

For a fixed $m\in\N$ define
\begin{equation*}
\Gamma^{(m)}=\{\gamma\in\Gamma:\gamma\equiv\text{Id (mod $m$)}\}.
\end{equation*}
We observe that $\Gamma$ acts on on the finite set $(\Z\slash
m\Z)^d.$ We denote by $\gamma\mapsto\bar{\gamma}$ the
corresponding homomorphism of $\Gamma$ into the semigroup
$\Lambda_{m,d}$ of maps of $(\Z\slash m\Z)^d$ into itself and we
write:
\begin{equation*}
\Gamma_m=\{\bar{\gamma}\in\Lambda_{m,d}:\gamma\in\Gamma\}.
\end{equation*}

The proof depends on the following
\begin{lemma}\label{lgamma}
Assume that $\Gamma$ is finitely generated and satisfies $(H_0),$
$(H_1)$ and $(H_2).$ Let $m$ be a prime number not dividing the
elements of the multiplicative semigroup
$\{\det\gamma:\gamma\in\Gamma\}.$ Then $\Gamma_m$ is a group of
permutations of $(\Z\slash m\Z)^d$ and the semigroup
$\Gamma^{(m)}$ satisfies $(H_0),$ $(H_1)$ and $(H_2).$
\end{lemma}
\begin{proof}
Here $\Z\slash m\Z$ is a finite field and for $\gamma\in\Gamma,$
$\bar{\gamma}$ is an endomorphism of the vector space $(\Z\slash
m\Z)^d.$ Then $\det\bar{\gamma}$ is the congruence class of
$\det\gamma$ in $\Z\slash m\Z.$ Since $m$ is a prime number not
dividing $\det\gamma,$ we conclude that
$\det\bar{\gamma}\not=\bar{0},$ hence $\bar{\gamma}$ belongs to
the group $GL(d,\Z\slash m\Z).$ Then $\Gamma_m$ is a semigroup
contained in the finite group $GL(d,\Z\slash m\Z);$ it follows
that $\Gamma_m$ is a group. We write
$\Gamma_m=\{\bar{a}_i:a_i\in\Gamma:i=1,\ldots,q\},$ and we observe
that the inverse of $\bar{a}_i$ is of the form
$\bar{a}_{i^\prime}$ with $a_{i^\prime}\in\Gamma$ and $1\leq
i^\prime\leq q.$ Since for every $\gamma\in\Gamma,$ we have
$\bar{\gamma}=\bar{a}_i$ for some $i,$ we get
$\bar{a}_{i^\prime}\bar{\gamma}=\text{Id},$
$a_{i^\prime}\gamma\in\Gamma^{(m)}.$

Assume condition $(H_1)$ is not satisfied by $\Gamma^{(m)};$ then
for some subspace $W\subset V,$ the orbit $\Gamma^{(m)}W$ is
finite, hence the set $\{a_{i^\prime}\gamma
W:\bar{a}_{i^\prime}\in\Gamma_m,\,\gamma\in\Gamma,\,\bar{a}_{i^\prime}\bar{\gamma}=\text{Id}\}$
is also finite. It follows that the set $\{\gamma
W:\gamma\in\Gamma\}$ is finite and this contradicts condition
$(H_1)$ for $\Gamma.$ Hence $\Gamma^{(m)}$ satisfies condition
$(H_1).$

Let $\gamma\in\Gamma$ be a proximal and expanding element element
of $\Gamma.$ Since the group $\Gamma_m$ is finite, there exists
$k\leq |\Gamma_m|$ such that $\bar{\gamma}^k=\text{Id},$ hence
$\gamma^k\in\Gamma^{(m)}.$ Clearly, $\gamma^k$ is proximal and
expanding. Then the result follows from the equivalence of a) and
c) in Proposition \ref{dodatkowa}.
\end{proof}
The following lemma will be used also. Its proof is analogous to
the classical case of one endomorphism of $\T^d$ (see for example
\cite{BM}). In this lemma, the torus $\T^d$ is endowed with its
normalized Haar measure, which is $\Gamma$-invariant.
\begin{lemma}\label{lemmaBM}
Assume $\Gamma\subset\m$ and $\Sigma\subset\T^d$ is measurable,
has positive measure and satisfies $\Gamma\Sigma\subset\Sigma.$
Then, if any character $\chi\not=\text{Id}$ has unbounded
$\Gamma^t$-orbit, then $\Sigma$ has measure $1,$ in particular
$\Gamma$ is ergodic on $\T^d.$
\end{lemma}
In order to prove Theorem \ref{main}, we can suppose $\Gamma$ to
be finitely generated. In fact, Proposition \ref{prop24} implies
that $\Gamma$ contains a finitely generated semigroup $\Gamma_1$
which satisfies $(H_1)$ and $(H_2),$ and we can add to $\Gamma_1$
an expanding element $\gamma$ from $\Gamma.$ Then the semigroup
generated by $\Gamma_1$ and $\gamma$ satisfies $(H_0),$ $(H_1)$
and $(H_2)$ in view of equivalence (c), (d) in Proposition
\ref{dodatkowa}.

Since $\Sigma$ is infinite and closed, it contains limit points.
We have two cases.

\textsf{Case 1.} Some limit point of $\Sigma$ is rational. So, let
$p\slash q$ be a limit point of $\Sigma.$ Then the set $q\Sigma$
is $\Gamma$-invariant and has $0$ as its limit point. Therefore,
by considerations in subsection \ref{subs1} we get that
$q\Sigma=\T^d.$ Hence, $\Sigma$ has positive Haar measure (greater
than $(1\slash q)^d$). Since $\Gamma$ satisfies $(H_0),$ $(H_1)$
we get, from Remark \ref{requiv} (iv) that $\Gamma^t$ satisfies
$(H_0),$ hence Lemma \ref{lemmaBM} allows us to conclude that
$\Sigma$ has measure $1.$ Since $\Sigma$ is closed, we have
$\Sigma=\T^d.$

\textsf{Case 2.} Every limit point of $\Sigma$ is irrational. Let
$\Sigma^{\text{ac}}$ be the set of limit points of $\Sigma.$ For
$m$ fixed and prime not dividing the elements of the finitely
generated semigroup $\{\det\gamma:\gamma\in\Gamma\},$ let
$\Sigma^{(m)}\subset\Sigma^{\text{ac}}$ be a minimal
$\Gamma^{(m)}$-invariant set. Since $\Sigma^{\text{ac}}$ consists
of irrational points, $\Sigma^{(m)}$ is infinite, hence $0$ is a
limit point of the closed $\Gamma^{(m)}$-invariant subset
$\Sigma^{(m)}-\Sigma^{(m)}.$ From Lemma \ref{lgamma} above and
consideration in subsection \ref{subs1} we get that
$\Sigma^{(m)}-\Sigma^{(m)}=\T^d.$ Therefore, for every
$r=(r_1,\ldots,r_d)\in\Z^d$ there are $x$ and $y$ in
$\Sigma^{(m)}$ such that
\begin{equation*}
x-y=(r_1\slash m,\ldots,r_d\slash m)=r\slash m.
\end{equation*}
Let $\Sigma^{(m)}_r$ be defined as follows:
\begin{equation*}
\Sigma^{(m)}_r=\{x\in\Sigma^{(m)}:\exists\,y\in\Sigma^{(m)}:x-y=r\slash
m)\}.
\end{equation*}
Clearly, $\Sigma^{(m)}_r$ is closed and nonempty. Since $r\slash
m$ is fixed by $\Gamma^{(m)}$ it follows that $\Sigma^{(m)}_r$ is
$\Gamma^{(m)}$-invariant. Thus, by minimality of $\Sigma^{(m)}$ we
get that $\Sigma^{(m)}=\Sigma^{(m)}_r.$ Therefore, for every
$m\in\N,$ $x\in\Sigma^{(m)},$ $r\in\Z^d$ we have
\begin{equation*}
x+r\slash m=y\in\Sigma^{(m)}.
\end{equation*}
Hence $\Sigma^{(m)}$ is invariant under translations in $\T^d$ by
$r\slash m,$ $r\in\Z^d.$ It follows that $\Sigma^{(m)}$ is
$1\slash m$-dense, hence $\Sigma^{\text{ac}}$ is $1\slash m$-dense
for every prime $m$ as above. We observe that the set of such
primes is infinite, thus $1\slash m$ can be chosen arbitrary
small. Since $\Sigma^{\text{ac}}$ is closed we have
$\Sigma^{\text{ac}}=\T^d,$ which contradicts the hypothesis. Thus,
only case 1 is possible, hence $\Sigma=\T^d.$ $\hfill\Box$


\begin{thebibliography}{1}

\bibitem{AMS}
H.~Abels, G.~A.~Margulis and G.~A.~So\u\i fer.
\newblock Semigroups containing proximal linear maps.
\newblock{\em  Israel J. Math.} 91(1-3):1--30, 1995.

\bibitem{BM}
M.~B.~Bekka and M.~Mayer.
\newblock Ergodic theory and topological dynamics of group actions on homogeneous
spaces. London Mathematical Society Lecture Note Series, 269, {\em
Cambridge University Press}, 2000.

\bibitem{Benoist}
Y.~Benoist.
\newblock Propri\'et\'es asymptotiques des groupes
lin\'eaires. II.
\newblock In {\em Analysis on homogeneous spaces and
representation theory of Lie groups, Okoyama - Kyoto 1997},
33--48, Adv. Stud. Pure Math., 26, Math. Soc. Japan, Tokyo, 2000.

\bibitem{B}
D.~Berend.
\newblock Multi-invariant sets on tori.
\newblock {\em Trans. Amer. Math. Soc.}, 280(2):509--532, 1983.

\bibitem{Bow}
R.~Bowen.
\newblock {\em Equilibrium states and the ergodic theory of Anosov
diffeomorphisms.}
\newblock Lecture Notes in Mathematics, Vol 470. Springer-Verlag, Berlin-New York, 1975.

\bibitem{CG00}
J.-P. Conze and Y.~Guivarc'h.
\newblock Limit sets of groups of linear transformations.
\newblock {\em Sankhy\=a Ser. A, Special issue on Ergodic Theory and Harmonic Analysis
in honour of M.G. Nadkarni}, 62(3):367--385, 2000.

\bibitem{CG01}
J.-P Conze and Y.~Guivarc'h.
\newblock Densit\'e d'orbites d'actions de groupes lin\'eaires et
  propri\'et\'es d'\'equidistribution de marches al\'eatoires.
\newblock In {\em Rigidity in dynamics and geometry (Cambridge, 2000)}, pages
  39--76. Springer, Berlin, 2002.

\bibitem{Dalbo}
F.~Dal'bo.
\newblock Topologie du feuilletage fortement stable.
\newblock {\em Ann. Inst. Fourier (Grenoble)}, 50(3):981--993, 2000.

\bibitem{DR}
S.~G.~Dani and S.~Raghavan.
\newblock Orbits of Euclidean frames
under discrete linear groups.
\newblock{\em Israel J. Math.}, 36(3-4):300--320, 1980.

\bibitem{SGP}
B. de Saporta, Y.~Guivarc'h and E.~Le Page.
\newblock On the multidimensional stochastic equation
$Y_{n+1}=A_nY_n+B_n.$
\newblock {\em C. R. Acad. Sci. Paris, Ser. I}, 339(7):499--502,
2004.

\bibitem{Eb}
P.~Eberlein.
\newblock Geodesic flows on negatively curved manifolds, II.
\newblock {\em Trans. Amer. Math. Soc.}, 178:57--82, 1973.

\bibitem{F}
H.~Furstenberg.
\newblock Disjointness in ergodic theory, minimal sets, and a problem in
{D}iophantine approximation.
\newblock {\em Math. Systems Theory}, 1:1--49, 1967.

\bibitem{F1}
H.~Furstenberg.
\newblock Boundary theory and stochastic processes on homogeneous
spaces. \newblock In {\em Harmonic Analysis on Homogeneous
Spaces}, Proc. Symp. Pure Math. 26:193--229, 1972, C.C. Moore
(ed.)

\bibitem{GG}
I.~Ya.~Goldsheid and Y.~Guivarc'h.
\newblock Zariski closure and the dimension of the Gaussian law of
the product of random matrices. I
\newblock {\em Probab. Theory Related Fields}, 105:109--142, 1996.

\bibitem{G}
Y.~Guivarc'h.
\newblock Produits de matrices al\'eatoires et
applications aux pro\-pri\-\'et\'es g\'e\-om\'e\-tri\-ques de
sous-groupes du groupe lin\'eaire.
\newblock {\em Ergodic Theory Dynam. Systems}, 10:483--512, 1990.

\bibitem{GLP}
Y.~Guivarc'h and E.~Le Page.
\newblock Simplicit\'e de spectres de Lyapunov et propri\'et\'e
d'isolation spectrale pour une famille d'op\'erateurs  de
transfert sur l'espace projectif.
\newblock In {\em Random Walks and Geometry,
Vienna Workshop, 2001}, 181--259, Ed. by Vadim Kaimanovich, Walter
de Gruyter, 2004.

\bibitem{GR}
Y.~Guivarc'h and A.~Raugi.
\newblock Products of random matrices: convergence theorems.
\newblock In {\em Random matrices and their applications (Brunswick, Maine,
  1984)}, volume~50 of {\em Contemp. Math.}, pages 31--54. Amer. Math. Soc.,
  Providence, RI, 1986.

\bibitem{GR1}
Y.~Guivarc'h and A.~Raugi.
\newblock Fronti\`ere de Furstenberg, propri\'et\'es de contraction et th\'eor\`emes de convergence
\newblock {\em Z. Wahrscheinlichkeitstheor. Verw. Geb.}, 69, 187--242, 1985.


\bibitem{GS}
Y.~Guivarc'h and A.~N. Starkov.
\newblock Orbits of linear group actions, random walk on homogeneous spaces,
  and toral automorphisms.
\newblock {\em Ergodic Theory Dynam. Systems}, 24(3):767--802, 2004.

\bibitem{HL}
G.~H. Hardy and J.~E. Littlewood.
\newblock Some problems of Diophantine approximation.
\newblock {\em Acta Math.}, 37:155--191, 1914.

\bibitem{Marinescu}
C.~T. Ionescu-Tulcea and G.~Marinescu.
\newblock Th\'e\-orie er\-go\-di\-que pour une classe d'op\'e\-ra\-te\-urs
non com\-pl\`e\-te\-ment continus
\newblock {\em Ann. of Math}, 52:140--147, 1950

\bibitem{K} H.~Kesten.
\newblock Random difference equations and renewal theory
for products of random matrices.
\newblock {\em Acta Math. 131 (1973), 207--248}.

\bibitem{L}
S.~P.~Lalley.
\newblock Renewal theorems in symbolic dynamics with
application to geodesic flows, non-Euclidean tessellations and
their fractal limits.
\newblock {\em Acta Math.}, 163(1-2):1--55, 1989.

\bibitem{Ma}
G.~A.~Margulis.
\newblock Problems and conjectures in rigidity theory.
\newblock In {\em Mathematics: frontiers and perspectives},
pages 161--174 of Amer. Math. Soc., Providence, RI, 2000

\bibitem{03a}
R.~Muchnik.
\newblock Orbits of Zariski dense semigroups of $SL(n,\mathbb Z)$.
\newblock {\em Ergodic Theory Dynam. Systems}, to appear

\bibitem{03b}
R.~Muchnik.
\newblock Semigroup actions on $\mathbb T^n$.
\newblock {\em Geometriae Dedicata}, to appear.

\bibitem{OV}
A.L.~Onischik and E.B.~Vinberg.
\newblock Lie groups and algebraic groups.
Springer Series in Soviet Mathematics.
\newblock{\em  Springer-Verlag, Berlin}, 1990

\bibitem{Ro}
M.~Rosenblatt.
\newblock Equicontinuous Markov operators
\newblock{\em Theory of Probability and its
Applications}, 2:180--197, 1964

\bibitem{S1}
A.~N.~Starkov.
\newblock Orbit closures of toral automorphism groups,
\newblock preprint, Moscow, 1999.

\end{thebibliography}
\end{document}